\documentclass[11pt]{article}
\usepackage{amssymb}

\usepackage[dvips]{graphicx}
\usepackage{amsmath}

\newtheorem{theorem}{Theorem}[subsection]

\newtheorem{corollary}[theorem]{Corollary}

\newtheorem{example}[theorem]{Example}

\newtheorem{lemma}[theorem]{Lemma}

\newtheorem{problem}[theorem]{Problem}
\newtheorem{proposition}[theorem]{Proposition}

\newenvironment{proof}[1][Proof]{\textbf{#1.} }{\ \rule{0.5em}{0.5em}}

\input{tcilatex}

\begin{document}

\title{Parabolic Kazhdan-Lusztig polynomials, plethysm and generalized
Hall-Littlewood functions for classical types}
\author{C\'{e}dric Lecouvey \\
Laboratoire de Math\'{e}matiques Pures et Appliqu\'{e}es Joseph Liouville\\
B.P. 699 62228 Calais Cedex\\
Cedric.Lecouvey@lmpa.univ-littoral.fr}
\date{}
\maketitle

\begin{abstract}
We use power sums plethysm operators to introduce $H$ functions which
interpolate between the Weyl characters and the Hall-Littlewood functions $%
Q^{\prime }$ corresponding to classical Lie groups. The coefficients of
these functions on the basis of Weyl characters are parabolic
Kazhdan-Lusztig polynomials and thus, by works of Kashiwara and Tanisaki,
are nonnegative.\ We prove that they can be regarded as quantizations of
branching coefficients obtained by restriction to certain subgroups of Levi
type. The $H$ functions associated to linear groups coincide with the
polynomials introduced by Lascoux Leclerc and Thibon in [J. Math. Phys 38
(1996), 1041-1068].
\end{abstract}

\section{Introduction}

Given $\mu $ a partition with at most $n$ parts, the Hall Littlewood
function $Q_{\mu }^{\prime }$ can be defined by 
\begin{equation*}
Q_{\mu }^{\prime }=\sum_{\lambda }K_{\lambda ,\mu }(q)s_{\lambda }
\end{equation*}
where the sum runs over the partitions of length at most $n$, $K_{\lambda
,\mu }(q)$ is the Lusztig $q$-analogue of weight multiplicity associated
with $(\lambda ,\mu )$ and $s_{\lambda }$ the Schur function indexed by $%
\lambda $, that is the Weyl character of the irreducible finite dimensional $%
GL_{n}$-module $V(\lambda )$. Since $K_{\lambda ,\mu }(1)$ is equal to the
dimension of the weight space $\mu $ in $V(\lambda )$, $Q_{\mu }^{\prime }$
can be regarded as a quantization of the homogeneous function $h_{\mu }.$ In 
\cite{LLT0}, Lascoux, Leclerc and Thibon have introduced a new family of
symmetric functions $H_{\mu }^{\ell }$ depending on a fixed nonnegative
integer $\ell $ which interpolate between the Schur functions $s_{\mu }$ and
the Hall-Littlewood functions $Q_{\mu }^{\prime }$. The polynomials $H_{\mu
}^{\ell }$ can be combinatorially described in terms of the spin statistic
on certain generalized Young tableaux called $\ell $-ribbon tableaux. These
ribbon tableaux naturally appear in the description of the action of the
power sum plethysm $\psi _{\ell }$ on symmetric functions.\ Recall that for
any symmetric function $f$, $\psi _{\ell }(f)$ is obtained by replacing in $%
f $ each variable $x_{i}$ by $x_{i}^{\ell }$. In particular $\psi _{\ell }$
multiplies the degrees by $\ell $. The space of symmetric functions is
endowed with an inner product $<\cdot ,\cdot >$ which makes the basis of
Schur functions orthonormal. Then $\varphi _{\ell }$, the adjoint operator
of $\psi _{\ell }$, divides the degrees by $\ell $. It is well known that $%
\varphi _{\ell }(s_{\mu })$ can be computed from the Jacobi-Trudi
determinantal identity.\ Namely, we have 
\begin{equation}
\varphi _{\ell }(s_{\mu })=0\text{ or }\varphi _{\ell }(s_{\mu
})=\varepsilon (\sigma _{0})s_{\mu ^{(0)}}\cdot \cdot \cdot s_{\mu ^{(\ell
-1)}}  \label{JT}
\end{equation}
where $\varepsilon (\sigma _{0})=\pm 1$ is the signature of a permutation $%
\sigma _{0}\in S_{n}$ and $(\mu ^{(0)},...,\mu ^{(\ell -1)})$ a $\ell $%
-tuple of partitions defined by $\ell $ and $\mu $. By expanding $\varphi
_{\ell }(s_{\mu })$ on the basis of Schur functions we obtain then 
\begin{equation}
\varphi _{\ell }(s_{\mu })=\varepsilon (\sigma _{0})\sum_{\lambda }c_{\mu
^{(0)},...,\mu ^{(\ell -1)}}^{\lambda }s_{\lambda }  \label{filA}
\end{equation}
where $c_{\mu ^{(0)},...,\mu ^{(\ell -1)}}^{\lambda }$ is the
Littlewood-Richardson coefficient giving the multiplicity of $V(\lambda )$
in the tensor product $V(\mu ^{(0)})\otimes \cdot \cdot \cdot \otimes V(\mu
^{(\ell -1})$. When $\ell =1,$ one has $\varphi _{\ell }(s_{\ell \mu
})=s_{\mu }$ and when $\ell >n$ one can prove that $\varphi _{\ell }(s_{\ell
\mu })=h_{\mu }$. Thus the functions $h_{\mu }^{(\ell )}=\varepsilon (\sigma
_{0})\varphi _{\ell }(s_{\ell \mu })$ interpolate between the functions $%
s_{\mu }$ and $h_{\mu }$ and have nonnegative coefficients on the basis of
Schur functions.

\noindent In \cite{LLT0}, the authors have interpreted the algebra of
symmetric functions as the bosonic Fock space representation of the quantum
affine Lie algebra $U_{q}(\widehat{sl_{n}})$.\ This permits to introduce a
natural quantization $\psi_{q,\ell}$ of the power sum plethysm $\psi_{\ell}$%
. Let $\varphi_{q,\ell}$ be the adjoint operator of $\psi_{q,\ell}$ with
respect to $<\cdot,\cdot>$. The function $H_{\mu}^{\ell}$ is then defined as
a simple renormalization of $\varphi_{q,\ell}(s_{\ell\mu})$. This gives an
identity of the form 
\begin{equation*}
H_{\mu}^{\ell}=\sum_{\lambda}c_{\mu^{(0)},...,\mu^{(\ell-1)}}^{\lambda
}(q)s_{\lambda}
\end{equation*}
where the polynomial $c_{\mu^{(0)},...,\mu^{(\ell-1)}}^{\lambda}(q)$ is a $q$%
-analogue of $c_{\mu^{(0)},...,\mu^{(\ell-1)}}^{\lambda}$.

\noindent Lusztig's $q$-analogues $K_{\lambda ,\mu }(q)$ are particular
affine Kazhdan-Lusztig polynomials. These polynomials arise in affine Hecke
algebra theory as the entries of the transition matrix between the natural
basis and a special basis defined by Lusztig.\ By replacing the affine Hecke
algebra $\widehat{H}$ by one of its parabolic modules $\widehat{H}\nu $ ($%
\nu $ being a weight of the affine root system under consideration), Deodhar
has introduced analogues of the Kazhdan-Lusztig polynomials.\ In \cite{LT},
it is shown that the family constituted by these parabolic Kazhdan-Lusztig
polynomials contains in particular the $q$-analogues $c_{\mu ^{(0)},...,\mu
^{(\ell -1)}}^{\lambda }(q)$.\ By a result of Kashiwara and Tanisaki \cite
{KT}, this implies notably that the coefficients of the polynomial $c_{\mu
^{(0)},...,\mu ^{(\ell -1)}}^{\lambda }(q)$ are nonnegative integers.

\bigskip

The aim of the paper is to introduce analogues of the polynomials $H_{\mu
}^{\ell }$ for the classical Lie groups $G=SO_{2n+1},Sp_{2n}$ and $SO_{2n}$
which interpolate between the Weyl characters and the Hall-Littlewood
functions associated with $G$. Write also $s_{\lambda }$ for the Weyl
character of the irreducible $G$-module $V(\lambda )$ of highest weight $%
\lambda $.\ We define the plethysm operator $\varphi _{\ell }$ and its dual $%
\psi _{\ell }$ on the $\mathbb{Z}$-algebra generated by these Weyl
characters.\ By a subgroup $L\subset G$ of Levi type, we mean a subgroup of $%
G$ isomorphic to the Levi subgroup of one of its parabolic subgroups.\ Given 
$\gamma $ a highest weight of $L$, we denote by $[V(\lambda ):V_{L}(\gamma
)] $ the multiplicity of the irreducible $L$-module $V_{L}(\gamma )$ of
highest weight $\gamma $ in the restriction of $V(\lambda )$ to $L$. Then,
provided that $\ell $ is odd when $G=Sp_{2n}$ or $SO_{2n}$, we establish for
any Weyl character $s_{\mu }$ such that $\varphi _{\ell }(s_{\mu })\neq 0$,
a formula of the type 
\begin{equation}
\varphi _{\ell }(s_{\mu })=\varepsilon (w_{0})\sum_{\lambda }[V(\lambda
):V_{L}\binom{\mu }{\ell }]s_{\lambda }  \label{filG}
\end{equation}
where $\varepsilon (w_{0})$ is the signature of an element $w_{0}\in W$ the
Weyl group of $G$, $L$ a subgroup of Levi type of $G$ and $\binom{\mu }{\ell 
}$ a dominant weight associated with $L$. The procedure which yields $%
w_{0},L $ and $\binom{\mu }{\ell }$ from $\ell $ and $\mu $ can be regarded
as an analogue of the algorithm computing the $\ell $-quotient of a
partition which implicitly appears in (\ref{JT}).\ The identity (\ref{filA})
can also be rewritten as in (\ref{filG}).\ Indeed, take $L=GL_{r_{0}}\times
\cdot \cdot \cdot \times GL_{r_{\ell -1}}$ where for any $k=1,...,\ell -1,$ $%
r_{k}$ is the length of $\mu ^{(k)}$.\ Then $\binom{\mu }{\ell }=(\mu
^{(0)},...,\mu ^{(\ell -1)})$ can be interpreted as a dominant weight for
the subgroup of Levi type $L$ of $GL_{n}$ and we have the duality $c_{\mu
^{(0)},...,\mu ^{(\ell -1)}}^{\lambda }=[V(\lambda ):V_{L}\binom{\mu }{\ell }%
]$.

\noindent The surprising constraint $\ell $ odd when $G=Sp_{2n}$ or $SO_{2n}$
appearing in (\ref{filG}) follows from the fact that the procedure giving $%
w_{0},L$ and $\binom{\mu }{\ell }$ mentioned above depends not only on the
Lie group under consideration but also on the parity of the integer $\ell $.
For $G=SO_{2n+1}$ the coefficients of $\varepsilon (w_{0})\varphi _{\ell
}(s_{\mu })$ on the basis of Weyl characters are always branching
coefficients corresponding to restriction to $L.$ For $G=Sp_{2n}$ or $%
SO_{2n} $ this is only true when $\ell $ is odd.\ Note that this difficulty
disappears in large rank, that is for $n\geq \ell \left| \mu \right| $ (but
see Section \ref{subsec-stabipleth}).

\noindent To define the functions $H_{\mu }^{\ell }$ in type $B,C$ or $D$,
we prove the equalities 
\begin{equation*}
\left| <\psi _{\ell }(s_{\lambda }),s_{\mu }>\right| =\left| <s_{\lambda
},\varphi _{\ell }(s_{\mu })>\right| =P_{\mu +\rho ,\ell \lambda +\rho
}^{-}(1)
\end{equation*}
which show that the coefficients of the expansion of $\psi _{\ell
}(s_{\lambda })$ on the basis of Weyl characters are, up to a sign,
parabolic Kazhdan-Lusztig polynomials specialized at $q=1.$ By using (\ref
{filG}) this gives, providing $\ell $ is odd for $G=Sp_{2n}$ or $SO_{2n}$%
\begin{equation*}
\lbrack V(\lambda ):V_{L}\binom{\mu }{\ell }]=P_{\mu +\rho ,\ell \lambda
+\rho }^{-}(1).
\end{equation*}
We then introduce the functions 
\begin{equation*}
G_{\mu }^{\ell }=\sum_{\lambda }[V(\lambda ):V_{L}\binom{\mu }{\ell }%
]_{q}s_{\lambda }
\end{equation*}
where $[V(\lambda ):V_{L}\binom{\mu }{\ell }]_{q}=P_{\mu +\rho ,\ell \lambda
+\rho }^{-}(q)$.\ This yields nonnegative $q$-analogues of the branching
coefficients $[V(\lambda ):V_{L}\binom{\mu }{\ell }]$.\ The functions $%
H_{\mu }^{\ell }$ are then defined by setting $H_{\mu }^{\ell }=G_{\ell \mu
}^{\ell }.$ We obtain the identities $H_{\mu }^{1}=s_{\mu }$ and $H_{\mu
}^{\ell }=Q_{\mu }^{\prime }$ when $\ell $ is sufficiently large.\ Thus the
functions $H_{\mu }^{\ell }$ interpolate between the Weyl characters and the
Hall-Littlewood functions associated with $G$.

\bigskip

The paper is organized as follows.\ In Section $2$ we recall the necessary
background on classical root systems, Weyl characters, subgroups of Levi
type and their corresponding branching coefficients. In Section $3,$ we
define the plethysm operators $\psi _{\ell }$ and their dual operators $%
\varphi _{\ell }.$ By abuse of notation, we also denote by $\varphi _{\ell }$
the linear operator on the group algebra $\mathbb{Z}[\mathbb{Z}^{n}]$ with
basis the formal exponentials $(e^{\beta })$ such that 
\begin{equation*}
\varphi _{\ell }(e^{\beta })=\left\{ 
\begin{array}{l}
e^{\beta /\ell }\text{ if }\beta \in (\ell \mathbb{Z})^{n} \\ 
0\text{ otherwise}
\end{array}
\right. \text{ for any }\beta \in \mathbb{Z}^{n}\text{.}
\end{equation*}
We then show how the determination of $\varphi _{\ell }(s_{\mu })$ can be
reduced to the computation of the polynomial 
\begin{equation*}
\varphi _{\ell }(e^{\mu }\prod_{\alpha \in R_{+}}(1-e^{\alpha }))
\end{equation*}
where $R_{+}$ is the set of positive roots corresponding to the Lie group $G$%
. This permits to establish formulas (\ref{filG}) providing $\ell $ is odd
when $G=Sp_{2n}$ or $SO_{2n}$. For completion we have also included the case 
$G=GL_{n}$ and shown why (\ref{filG}) cannot hold when $\ell $ is even and $%
G=Sp_{2n}$ or $SO_{2n}$. To make the paper self-contained, we have
summarized in Section $4$ some necessary results on affine Hecke algebras
and parabolic Kazhdan-Lusztig polynomials. Section $5$ is devoted to the
definition of the polynomials $G_{\mu }^{\ell }$ and $H_{\mu }^{\ell }$ and
to their links with the Weyl characters and the Hall-Littlewood functions.\
Finally we briefly discuss in Section $6$ the problem of defining
nonnegative $q$-analogues of tensor product multiplicities when $G\neq
GL_{n} $. We add also a few remarks concerning the exceptional root systems.

\bigskip

\noindent \textbf{Acknowledgments: }The author wants to thank B. Leclerc for
very helpful and stimulating discussions on the results of \cite{LLT0} and 
\cite{LT}.

\section{Background}

\subsection{Classical root systems}

In the sequel $G$ is one of the complex Lie groups $GL_{n},Sp_{2n},SO_{2n+1}$
or $SO_{2n}$ and $\frak{g}$ its Lie algebra.\ We follow the convention of 
\cite{KT} to realize $G$ as a subgroup of $GL_{N}$ and $\frak{g}$\ as a
subalgebra of $\frak{gl}_{N}$ where 
\begin{equation*}
N=\left\{ 
\begin{tabular}{l}
$n$ when $G=GL_{n}$ \\ 
$2n$ when $G=Sp_{2n}$ or $SO_{2n}$ \\ 
$2n+1$ when $G=SO_{2n+1}$%
\end{tabular}
\right. .
\end{equation*}
With this convention the maximal torus $T$ of $G$ and the Cartan subalgebra $%
\frak{h}$ of $\frak{g}$ coincide respectively with the subgroup and the
subalgebra of diagonal matrices of $G$ and $\frak{g}$. Similarly the Borel
subgroup $B$ of $G$ and the Borel subalgebra $\frak{b}$ of $\frak{g}$
coincide respectively with the subgroup and subalgebra of upper triangular
matrices of $G$ and $\frak{g}$.

\noindent Let $d_{N}$ be the linear subspace of $\frak{gl}_{N}$ consisting
of the diagonal matrices.\ For any $i\in I_{n}=\{1,...,n\},$ write $%
\varepsilon _{i}$ for the linear map $\varepsilon _{i}:d_{N}\rightarrow 
\mathbb{C}$ such that $\varepsilon _{i}(D)=\delta _{n-i+1}$ for any diagonal
matrix $D$ whose $(i,i)$-coefficient is $\delta _{i}.$ Then $(\varepsilon
_{1},...,\varepsilon _{n})$ is an orthonormal basis of the Euclidean space $%
\frak{h}_{\mathbb{R}}^{\ast }$ (the real part of $\frak{h}^{\ast }).$ Let $%
(\cdot ,\cdot )$ be the corresponding nondegenerate symmetric bilinear form
defined on $\frak{h}_{\mathbb{R}}^{\ast }$.\ Write $R$ for the root system
associated with $G.\;$For any $\alpha \in R$ we set $\alpha ^{\vee }=\frac{%
\alpha }{(\alpha ,\alpha )}$. The Lie algebra $\frak{g}$ admits the diagonal
decomposition $\frak{g}=\frak{h}\oplus \bigoplus_{\alpha \in R}\frak{g}%
_{\alpha }.\;$We take for the set of positive roots: 
\begin{equation*}
\left\{ 
\begin{tabular}{l}
$R^{+}=\{\varepsilon _{j}-\varepsilon _{i}\text{ with }1\leq i<j\leq n\}%
\text{ for the root system }A_{n-1}$ \\ 
$R^{+}=\{\varepsilon _{j}-\varepsilon _{i},\varepsilon _{j}+\varepsilon _{i}%
\text{ with }1\leq i<j\leq n\}\cup \{\varepsilon _{i}\text{ with }1\leq
i\leq n\}\text{ for the root system }B_{n}$ \\ 
$R^{+}=\{\varepsilon _{j}-\varepsilon _{i},\varepsilon _{j}+\varepsilon _{i}%
\text{ with }1\leq i<j\leq n\}\cup \{2\varepsilon _{i}\text{ with }1\leq
i\leq n\}\text{ for the root system }C_{n}$ \\ 
$R^{+}=\{\varepsilon _{j}-\varepsilon _{i},\varepsilon _{j}+\varepsilon _{i}%
\text{ with }1\leq i<j\leq n\}\text{ for the root system }D_{n}$%
\end{tabular}
\right. .
\end{equation*}
Let $\rho $ be the half sum of positive roots.\ Set $J_{n}=\{\overline{n}%
<\cdot \cdot \cdot <\overline{1}<1<\cdot \cdot \cdot <n\}$ where, for each
integer $i=1,...,n$, we have written $\overline{i}$ for the negative integer 
$-i$.\ For any $x\in J_{n}$ we have $\overline{\overline{x}}=x$ and we set $%
\left| x\right| =x$ if $x>0,$ $\left| x\right| =\overline{x}$ otherwise.
Given a subset $U\subset J_{n},$ we define $\left| U\right| =\{\left|
x\right| \mid x\in U\}$ and $\overline{U}=\{\overline{x}\mid x\in U\}$.

\noindent The Weyl group of $GL_{n}$ is the symmetric group $S_{n}$ and for $%
G=SO_{2n+1},Sp_{2n}$ or $SO_{2n},$ the Weyl group $W$ of the Lie group $G$
is the subgroup of the permutation group of $J_{n}$ generated by the
permutations 
\begin{equation*}
\left\{ 
\begin{tabular}{l}
$s_{i}=(i,i+1)(\overline{i},\overline{i+1}),$ $i=1,...,n-1$ and $s_{n}=(n,%
\overline{n})$ $\text{for the root systems }B_{n}$ and $C_{n}$ \\ 
$s_{i}=(i,i+1)(\overline{i},\overline{i+1}),$ $i=1,...,n-1$ and $%
s_{n}^{\prime }=(n,\overline{n-1})(n-1,\overline{n})$ $\text{for the root
system }D_{n}$%
\end{tabular}
\right.
\end{equation*}
where for $a\neq b$ $(a,b)$ is the simple transposition which switches $a$
and $b.$ For types $B_{n}$ and $C_{n}$, $W$ is the group of signed
permutations.\ It is the subgroup of the permutation group of $J_{n}$
consisting of the permutations $w$ such that $w(\overline{i})=\overline{w(i)}
$.\ For type $D_{n},$ the elements of $W$ verify the additional constraint $%
\mathrm{card}\{i\in I_{n}\mid w(i)<0\}\in 2\mathbb{N}$. We identify the
subgroup of $W$ generated by $s_{i}=(i,i+1)(\overline{i},\overline{i+1}),$ $%
i=1,...,n-1$ with the symmetric group $S_{n}.$ The signature $\varepsilon $
of $w\in W$ is defined by $\varepsilon (w)=(-1)^{l(w)}$ where $l$ is the
length function corresponding to the above sets of generators. Consider the
increasing sequence $K=(\overline{i}_{p},...,\overline{i}%
_{1},i_{1},...,i_{p})\subset J_{n}$. For $X=B,D$ set 
\begin{equation*}
W_{X,K}=\{w\in W\text{ of type }X_{n}\mid w(x)=x\text{ for any }x\notin K\}.
\end{equation*}
Then, $W_{X,K}$ is isomorphic to the Weyl group of type $X_{p}.$ Let $%
\varepsilon _{X,K}$ be the corresponding signature.

\begin{lemma}
\label{lem_sign}Consider $X=B,D$ and $w\in W_{X,K}$.\ Then we have $%
\varepsilon _{X,K}(w)=\varepsilon (w)$.
\end{lemma}

\begin{proof}
Suppose $X=B$.\ The generators of the Weyl group $W_{X,K}$ are the $%
t_{k}=(i_{k},i_{k+1})(\overline{i}_{k},\overline{i}_{k+1}),$ $k=1,...,p-1$
and $s_{n}=(\overline{i}_{p},i_{p})$.\ One verifies easily that, considered
as elements of $W$, they have an odd length. We proceed similarly when $X=D$.
\end{proof}

\bigskip

\noindent The action of $w\in W$ on $\beta =(\beta _{1},...,\beta _{n})\in 
\frak{h}_{\mathbb{R}}^{\ast }$ is defined by 
\begin{equation}
w\cdot (\beta _{1},...,\beta _{n})=(\beta _{1}^{w^{-1}},...,\beta
_{n}^{w^{-1}})  \label{actionW}
\end{equation}
where $\beta _{i}^{w}=\beta _{w(i)}$ if $w(i)\in I_{n}$ and $\beta
_{i}^{w}=-\beta _{w(\overline{i})}$ otherwise. The dot action of $W$ on $%
\beta =(\beta _{1},...,\beta _{n})\in \frak{h}_{\mathbb{R}}^{\ast }$ is
defined by 
\begin{equation}
w\circ \beta =w\cdot (\beta +\rho )-\rho .  \label{dotaction}
\end{equation}
\ The fundamental weights of $\frak{g}$ belong to $\left( \frac{\mathbb{Z}}{2%
}\right) ^{n}$.\ More precisely we have $\omega _{i}=(0^{i},1^{i})\in 
\mathbb{N}^{n}$ for $i<n-1$ and also $i=n-1$ for $\frak{g}\neq \frak{so}%
_{2n} $, $\omega _{n}^{C_{n}}=(1^{n}),$ $\omega _{n}^{B_{n}}=\omega
_{n}^{D_{n}}=(\frac{1}{2}^{n})$ and $\omega _{n-1}^{D_{n}}=(-\frac{1}{2},%
\frac{1}{2}^{n-1}).$ The weight lattice $P$ of $\frak{g}$ can be considered
as the $\mathbb{Z}$-sublattice of $\left( \frac{\mathbb{Z}}{2}\right) ^{n}$
generated by the $\omega _{i},$ $i\in I.$ For any $\beta =(\beta
_{1},...,\beta _{n})\in P,$ we set $\left| \beta \right| =\beta _{1}+\cdot
\cdot \cdot +\beta _{n}.\;$Write $P^{+}$ for the cone of dominant weights of 
$G.\;$With our convention, a partition of length $m$ is a weakly \textit{%
increasing} sequence of $m$ nonnegative integers.\ Denote by $\mathcal{P}%
_{n} $ the set of partitions with at most $n$ parts. Each partition $\lambda
=(\lambda _{1},...,\lambda _{m})\in \mathcal{P}_{n}$ will be identified with
the dominant weight $\sum_{i=1}^{m}\lambda _{i}\varepsilon _{i}.\;$Then the
irreducible finite dimensional polynomial representations of $G$ are
parametrized by the partitions of $\mathcal{P}_{n}$.\ For any $\lambda \in 
\mathcal{P}_{n},$ denote by $V(\lambda )$ the irreducible finite dimensional
representation of $G$ of highest weight $\lambda .$ We will also need the
irreducible rational representations of $GL_{n}$.\ They are indexed by the $%
n $-tuples 
\begin{equation}
(\gamma ^{-},\gamma ^{+})=(-\gamma _{q}^{-},...,-\gamma _{1}^{-},\gamma
_{1}^{+},\gamma _{2}^{+},...,\gamma _{p}^{+})  \label{jamma+-}
\end{equation}
where $\gamma ^{+}=(\gamma _{1}^{+},\gamma _{2}^{+},...,\gamma _{p}^{+})$
and $\gamma ^{-}=(\gamma _{1}^{-},...,\gamma _{q}^{-})$ are partitions of
length $p$ and $q$ such that $p+q=n.$ Write $\widetilde{\mathcal{P}}_{n}$
for the set of such $n$-tuples and denote also by $V(\gamma )$ the
irreducible rational representation of $GL_{n}$ of highest weight $\gamma
=(\gamma ^{-},\gamma ^{+})\in \widetilde{\mathcal{P}}_{n}.$

\bigskip

\noindent In the sequel, our computations will also make appear root
subsystems of the root systems $R$ described above.\ Suppose that $G$ is of
type $X_{n}$ with $X_{n}=A_{n-1},B_{n},C_{n}$ or $D_{n}$.\ Let $%
I=(i_{1},...,i_{r})$ be an increasing sequence of integers belonging to $%
I_{n}$, that is $i_{k}\in I_{n}$ for any $k=1,..,r$ and $i_{1}<\cdot \cdot
\cdot <i_{r}$.\ Then 
\begin{equation*}
R_{I}=\{\alpha \in R\cap \oplus _{i\in I}\mathbb{Z}\varepsilon _{i}\}
\end{equation*}
is a root subsystem of $R$ of type $X_{r}$.\ Write $R_{I}^{+}$ for the set
of positive roots in $R_{I}.\;$Then we have $R_{I}^{+}=R_{I}\cap R^{+}.$ The
dominant weights associated with $R_{I}$ have the form $\lambda =(\lambda
_{1},...,\lambda _{n})$ where $\lambda _{i}\neq 0$ only if $i\in I$ and $%
\lambda ^{(I)}=(\lambda _{i_{1}},...,\lambda _{i_{r}})\in \mathcal{P}_{r}$.
We slightly abuse the notation by identifying $\lambda $ with $\lambda
^{(I)}.$

\noindent Consider an increasing sequence $X=(x_{1},...,x_{r})$ of integers
belonging to $J_{n}$ such that $\left| x_{k}\right| =\left|
x_{k^{\prime}}\right| $ if and only if $k=k^{\prime}$.\ For any integer $%
i=1,...,n,$ set $\varepsilon_{\overline{i}}=-\varepsilon_{i}.$ Then 
\begin{equation*}
R_{A,X}=\{\pm(\varepsilon_{x_{j}}-\varepsilon_{x_{i}})\mid1\leq i<j\leq r\}
\end{equation*}
is a root subsystem of $R$ of type $A_{r-1}.$ To see this, consider the
linear map $\theta_{X}:\mathbb{Z}^{r}\rightarrow\mathbb{Z}^{n}$ such that $%
\theta _{X}(\varepsilon_{i})=\varepsilon_{x_{i}}.$ The map $\theta$ is
injective and preserves the scalar product in $\mathbb{Z}^{r}$ and $\mathbb{Z%
}^{n}.$ Moreover the root system $\{\pm(\varepsilon_{j}-\varepsilon_{i})%
\mid1\leq i<j\leq r\}\subset\mathbb{Z}^{r}$ of type $A_{r}$ is sent on $%
R_{A,X}$ by $\theta_{X}.$ The set of positive roots in $R_{A,X}$ is equal to 
$R_{A,X}^{+}=R_{A,X}\cap R^{+}.$ Denote by $s\in\{1,...,r\}$ the maximal
integer such that $x_{s}<0$.\ We associate to $X,$ the increasing sequence
of indices $I\subset I_{n}$ defined by 
\begin{equation}
I=(\overline{x}_{s},...,\overline{x}_{1},x_{s+1},...,x_{r}).  \label{defIX}
\end{equation}
It will be useful to consider the weights corresponding to $R_{A,X}$ as the $%
r$-tuples $\beta=(\beta_{x_{1}},...,\beta_{x_{r}})$ with coordinates indexed
by $X.\;$The coordinates $(\beta_{1}^{\prime},...,\beta_{n}^{\prime})$ of $%
\beta$ on the initial basis $(\varepsilon_{1},...,\varepsilon_{n})$ are such
that $\beta_{i}^{\prime}=\beta_{x_{a}}$ if $i=x_{a}\in X,$ $%
\beta_{i}^{\prime }=-\beta_{x_{a}}$ if $\overline{i}=x_{a}\in X$ and $%
\beta_{i}^{\prime}=0$ otherwise. With this convention the dominant weights
for $R_{A,X}$ have the form 
\begin{equation}
\lambda^{(X)}=(\lambda_{x_{1}},...,\lambda_{x_{r}})\in\widetilde{\mathcal{P}}%
_{r}.  \label{conv}
\end{equation}
This simply means that we have chosen to expand the weights of $R_{A,X}$ on
the basis $\{\varepsilon_{x}\mid x\in X\}$ rather than on the basis $%
\{\varepsilon_{i}\mid i\in I\}$ to preserve the identification of the
dominant weights with the nondecreasing $r$-tuples of integers.

\begin{example}
\label{exam1}Take $G=Sp_{10}$.

\begin{itemize}
\item  For $I=(2,4,5)$ we have 
\begin{equation*}
R_{I}^{+}=\{\varepsilon _{5}\pm \varepsilon _{4},\varepsilon _{5}\pm
\varepsilon _{2},\varepsilon _{4}\pm \varepsilon _{2},2\varepsilon
_{2},2\varepsilon _{4},2\varepsilon _{5}\}
\end{equation*}
which is the set of positive roots of a root system of type $C_{3}.$ The
weight $\lambda =(1,2,2)$ is dominant for $G_{I}.\;$Considered as a weight
of $Sp_{10},$ we have $\lambda =(0,1,0,2,2).$

\item  For $X=(\overline{5},\overline{2},1,4)$ we have 
\begin{equation*}
R_{A,X}^{+}=\{\varepsilon _{4}-\varepsilon _{1},\varepsilon _{5}-\varepsilon
_{2},\varepsilon _{1}+\varepsilon _{2},\varepsilon _{1}+\varepsilon
_{5},\varepsilon _{4}+\varepsilon _{2},\varepsilon _{4}+\varepsilon _{5}\}
\end{equation*}
which is the set of positive roots of a root system of type $A_{3}.$ The
weight $\gamma =(-3,-1,4,5)$ is dominant for $G_{X}.\;$Considered as a
weight of $Sp_{10},$ we have $\gamma =(4,1,0,5,3).$
\end{itemize}
\end{example}

\subsection{Subgroups of Levi type\label{subsec}}

Suppose $G$ is a classical Lie group and consider $R$ the corresponding root
system. We shall need Lie subgroups of $G$ associated with particular
sub-root systems of $R.$ Each of these subgroups will be of Levy type, that
is, will be isomorphic to the Levi subgroup of one of the parabolic
subgroups of $G$.

\noindent Consider $p\geq 1$ an integer. Let $%
I^{(0)}=(i_{1}^{(0)},...,i_{r_{0}}^{(0)})$ be an increasing sequence of
integers in $I_{n}$.\ For $k=1,...,p,$ consider increasing sequences $%
X^{(k)}=(x_{1}^{(k)},...,x_{r_{k}}^{(k)})\subset J_{n}$ such that $\mathrm{%
card}(X^{(k)})=r_{k}.$ Let $s_{k}$ be maximal in $\{1,...,r_{k}\}$ such that 
$x_{s_{k}}^{(k)}<0$. Set 
\begin{equation}
I^{(k)}=(\overline{x}_{s_{k}}^{(k)},...,\overline{x}%
_{1}^{(k)},x_{s_{k}+1}^{(k)},...,x_{r_{k}}^{(k)})\subset I_{n}.  \label{f_Ik}
\end{equation}
We suppose that the sets $I^{(k)},$ $k=0,...,p$ are pairwise disjoint and
verify $\cup _{k=0}^{p}I^{(k)}=I_{n}.$ Set $\mathcal{I}%
=\{I^{(0)},X^{(1)},...,X^{(p)}\}$ and 
\begin{equation*}
R_{\mathcal{I}}=R_{I^{(0)}}\cup \bigcup_{k=1}^{p}R_{A,X^{(k)}}
\end{equation*}
Then $\frak{g}_{\mathcal{I}}=\frak{h}\oplus \bigoplus_{\alpha \in R_{%
\mathcal{I}}}\frak{g}_{\alpha }$ is a Lie subalgebra of $\frak{g}$.\ Its
corresponding Lie group $G_{\mathcal{I}}$ is a subgroup of $G$ of Levi type
and we have 
\begin{equation*}
G_{\mathcal{I}}\simeq \left\{ 
\begin{array}{l}
GL_{r_{0}}\times GL_{r_{1}}\times \cdot \cdot \cdot \times GL_{r_{p}}\text{
for }G=GL_{n} \\ 
SO_{2r_{0+1}}\times GL_{r_{1}}\times \cdot \cdot \cdot \times GL_{r_{p}}%
\text{ for }G=SO_{2n+1} \\ 
Sp_{2r_{0}}\times GL_{r_{1}}\times \cdot \cdot \cdot \times GL_{r_{p}}\text{
for }G=Sp_{2n} \\ 
SO_{2r_{0}}\times GL_{r_{1}}\times \cdot \cdot \cdot \times GL_{r_{p}}\text{
for }G=SO_{2n}
\end{array}
\right. .
\end{equation*}
The root system associated with $G_{\mathcal{I}}$ is $R_{\mathcal{I}}$.\
Denote by $P_{\mathcal{I}}^{+}$ its cone of dominant weights. The weight
lattice of $G_{\mathcal{I}}$ coincides with that\ of $G$ since the Lie
algebras $\frak{g}_{\mathcal{I}}$ and $\frak{g}$ have the same Cartan
subalgebra. The elements of $P_{\mathcal{I}}^{+}$ are the $(p+1)$-tuples $%
\lambda =(\lambda ^{(0)},\lambda ^{(1)},...,\lambda ^{(p)})$ where $\lambda
^{(0)}=(\lambda _{i}\mid i\in I^{(0)})$ is a dominant weight of $%
R_{G,I^{(0)}}$ and for any $k=1,...,p,$ $\lambda ^{(k)}=(\lambda _{i}\mid
i\in X^{(k)})$ is a dominant weight of $R_{G,X^{(k)}}$.\ For any $\lambda
\in P_{\mathcal{I}}^{+}$, we denote by $V_{\mathcal{I}}(\lambda )$ the
irreducible finite dimensional $G_{\mathcal{I}}$-module of highest weight $%
\lambda $. Each weight $\beta =(\beta ^{(0)},\beta ^{(1)},...,\beta
^{(p)})\in P_{\mathcal{I}}$ can be considered as a weight $\beta =(\beta
_{1}^{\prime },...,\beta _{n}^{\prime })$ of $P$. With the convention (\ref
{conv}) we have then $\beta _{i}^{\prime }=\beta _{i_{a}^{(0)}}$ if $%
i=i_{a}^{(0)}\in I^{(0)}$ and for any $k=1,...,p$, $\beta _{i}^{\prime
}=\beta _{i_{a}^{(k)}}$ if $i=i_{a}^{(k)}\in X^{(k)},$ $\beta _{i}^{\prime
}=-\beta _{i_{a}^{(k)}}$ if $\overline{i}=i_{a}^{(k)}\in X^{(k)}$. In the
sequel we identify the two expressions 
\begin{equation}
\beta =(\beta ^{(0)},\beta ^{(1)},...,\beta ^{(p)})\text{ and }\beta =(\beta
_{1}^{\prime },...,\beta _{n}^{\prime })  \label{iden}
\end{equation}
of the weights of $P_{\mathcal{I}}$.

\subsection{Weyl characters and dual bases}

We refer the reader to \cite{mac} and \cite{Ram} for a detailed exposition
of the results used in this paragraph. We use as a basis of the group
algebra $\mathbb{Z}[\mathbb{Z}^{n}],$ the formal exponentials $(e^{\beta
})_{\beta \in \mathbb{Z}^{n}}$ satisfying the relations $e^{\beta
_{1}}e^{\beta _{2}}=e^{\beta _{1}+\beta _{2}}.$ We furthermore introduce $n$
independent indeterminates $x_{1},...,x_{n}$ in order to identify $\mathbb{Z}%
[\mathbb{Z}^{n}]$ with the ring of polynomials $\mathbb{Z}%
[x_{1},...,x_{n},x_{1}^{-1},...,x_{n}^{-1}]$ by writing $e^{\beta
}=x_{1}^{\beta _{1}}\cdot \cdot \cdot x_{n}^{\beta _{n}}=x^{\beta }$ for any 
$\beta =(\beta _{1},...,\beta _{n})\in \mathbb{Z}^{n}.$ Define the action of
the Weyl group $W$ on $\mathbb{Z}[\mathbb{Z}^{n}]$ by $w\cdot x^{\beta
}=x^{w(\beta )}.$ The Weyl character $s_{\beta }$ is defined by 
\begin{equation*}
s_{\beta }=\dfrac{a_{\beta +\rho }}{a_{\rho }}\text{ where }a_{\beta
}=\sum_{w\in W}\varepsilon (\sigma )(w\cdot x^{\beta })
\end{equation*}
For any $\beta \in \mathbb{Z}^{n}$ we have 
\begin{equation}
s_{\beta }=\left\{ 
\begin{tabular}{l}
$\varepsilon (w)s_{\lambda }$ if there exists $w\in W$ and $\lambda \in 
\mathcal{P}_{n}$ such that $\lambda =w\circ \beta $ \\ 
$0$ otherwise
\end{tabular}
\right. .  \label{SL}
\end{equation}
Let $A$ be the $\mathbb{Z}$-algebra generated by the characters $s_{\lambda
},\lambda \in \mathcal{P}_{n}$.\ For any $\beta \in \mathbb{Z}^{n}$, denote
by $W_{\beta }$ the stabilizer of $\beta $ under the action of the Weyl
group $W$ and by $W^{\beta }$ a set of representatives in $W/W_{\beta }$
with minimal length.\ Then the functions 
\begin{equation*}
m_{\beta }=\sum_{w\in W^{\beta }}w\cdot x^{\beta }
\end{equation*}
belong to $A$. Moreover $\{m_{\lambda }\mid \lambda \in \mathcal{P}_{n}\}$
is a basis of $A.\;$We have the decomposition 
\begin{equation}
s_{\lambda }=\sum_{\mu \in \mathcal{P}_{n}}K_{\lambda ,\mu }m_{\mu }
\label{c_s}
\end{equation}
where $K_{\lambda ,\mu }$ is equal to the dimension of the weight space $\mu 
$ in the irreducible representation $V(\lambda )$. There exists an inner
product $<\cdot ,\cdot >$ on $A$ which makes the characters $s_{\lambda }$
orthonormal.\ We denote by $\{h_{\mu }\mid \mu \in \mathcal{P}_{n}\}$ the
dual basis of $\{m_{\lambda }\mid \lambda \in \mathcal{P}_{n}\}$ with
respect to $<\cdot ,\cdot >$. The homogeneous functions $h_{\mu }$ are given
in terms of the Weyl characters by the decomposition 
\begin{equation}
h_{\mu }=\sum_{\lambda \in \mathcal{P}_{n}}K_{\lambda ,\mu }s_{\lambda }.
\label{dec_h}
\end{equation}
This decomposition is infinite in general when $G\neq GL_{n}.$ Nevertheless,
by embedding $A$ in the ring $\widehat{A}$ of universal characters defined
by Koike and Terada \cite{KT}, it makes sense to consider formal series in
the characters $s_{\lambda },\lambda \in \mathcal{P}_{n}$. Note that the
function $h_{\mu }$ is not the character of the representation $V(\mu
_{1}\omega _{1})\otimes \cdot \cdot \cdot \otimes V(\mu _{n}\omega _{1})$
when $G\neq GL_{n}$, . For any $\beta \in \mathbb{Z}^{n}$ we define the
function $h_{\beta }$ by 
\begin{equation}
h_{\beta }=h_{\mu }  \label{SLH}
\end{equation}
where $\mu $ is the unique dominant weight contained in the orbit $W\cdot
\beta $.

\subsection{Jacobi-Trudi identities\label{subsec_jt}}

\noindent Denote by $\mathcal{L}_{n}=\mathbb{K[}[x^{\beta }]]$ the vector
space of formal series in the monomials $x^{\beta }$ with $\beta \in \mathbb{%
Z}$.\ We identify the ring of polynomials $\mathcal{F}_{n}=\mathbb{K[}%
x^{\beta }]$ with the subspace of $\mathcal{L}_{n}$ containing the finite
formal series.\ The vector space $\mathcal{L}_{n}$ is not a ring since $%
\beta \in \mathbb{Z}$. More precisely, the product $F_{1}\cdot \cdot \cdot
F_{r}$ of the formal series $F_{i}=\sum_{\beta \in E_{i}}x^{\beta ^{(i)}}$ $%
i=1,...,r$ is defined if and only if, for any $\gamma \in \mathbb{Z}^{n},$
the number $N_{\gamma }$ of decompositions $\gamma =\beta ^{(1)}+\cdot \cdot
\cdot +\beta ^{(r)}$ such that $\beta ^{(i)}\in E_{i}$ is finite and in this
case we have 
\begin{equation*}
F_{1}\cdot \cdot \cdot F_{r}=\sum_{\gamma \in \mathbb{Z}^{n}}N_{\gamma
}x^{\gamma }.
\end{equation*}
In particular the product $P\cdot F$ with $P\in \mathcal{P}_{n}$ and $F\in 
\mathcal{L}_{n}$ is well defined.

\noindent Set 
\begin{equation*}
\nabla=\prod_{\alpha\in R_{+}}\frac{1}{(1-x^{\alpha})}\text{ and }\Delta
=\prod_{\alpha\in R_{+}}(1-x^{\alpha}).
\end{equation*}
Then $\Delta\in\mathcal{F}_{n}$ and $\nabla\in\mathcal{L}_{n}$. We define
two linear maps 
\begin{equation*}
\mathrm{S}:\left\{ 
\begin{tabular}{c}
$\mathcal{L}_{n}\rightarrow\widehat{A}$ \\ 
$x^{\beta}\mapsto s_{\beta}$%
\end{tabular}
\right. \text{ and }\mathrm{H}:\left\{ 
\begin{tabular}{c}
$\mathcal{L}_{n}\rightarrow\widehat{A}$ \\ 
$x^{\beta}\mapsto h_{\beta}$%
\end{tabular}
\right. \text{.}
\end{equation*}
From Theorem 2.14 of \cite{Ram} we obtain

\begin{proposition}
\label{prop_JT}For any $\beta \in \mathbb{Z}^{n},$ $s_{\beta }=\sum_{w\in
W}\varepsilon (w)h_{\beta +\rho -w\cdot \rho }.$
\end{proposition}

\noindent By using the identity 
\begin{equation}
\Delta =x^{\rho }\sum_{w\in W}\varepsilon (w)x^{-w\cdot \rho }
\label{iden_delta}
\end{equation}
the previous proposition is equivalent to the following identity: 
\begin{equation}
\mathrm{S}(x^{\beta })=\mathrm{H}(\Delta \times x^{\beta }).  \label{S_H}
\end{equation}

\begin{proposition}
\label{propHS}For any $\beta \in \mathbb{Z}^{n}$ we have $\mathrm{H}%
(x^{\beta })=\mathrm{S}(\nabla \times x^{\beta }).$
\end{proposition}

\begin{proof}
Denote by $\chi _{\Delta }$ and $\chi _{\nabla }$ the linear maps defined on 
$\mathcal{L}_{n}$ by setting $\chi _{\Delta }(x^{\beta })=\Delta \times
x^{\beta }$ and $\chi _{\nabla }(x^{\beta })=\nabla \times x^{\beta }$
respectively. By (\ref{S_H}) we have $\mathrm{S}=\mathrm{H}\circ \chi
_{\Delta }.$ Moreover for any $\beta \in \mathbb{Z}^{n}$, $\chi _{\Delta
}\circ \chi _{\nabla }(x^{\beta })=x^{\beta }.$ This gives $\mathrm{S}%
(\nabla \times x^{\beta })=\mathrm{S}\circ \chi _{\nabla }(x^{\beta })=%
\mathrm{H}\circ \chi _{\Delta }\circ \chi _{\nabla }(x^{\beta })=\mathrm{H}%
(x^{\beta }).$
\end{proof}

\subsection{Branching coefficients for the restriction to subgroups of Levi
type\label{subsec_BC}}

\noindent Consider $\mathcal{I}=\{I_{0},X_{1},...,X_{p}\}$ as in \ref{subsec}%
. The set $\mathcal{I}$ characterizes a subgroup $G_{\mathcal{I}}\subset G$
of Levi type. Set 
\begin{equation*}
\Delta _{\mathcal{I}}=\prod_{\alpha \in R_{\mathcal{I}}^{+}}(1-x^{\alpha })%
\text{ and }\nabla _{\mathcal{I}}=\prod_{\alpha \in R^{+}-R_{\mathcal{I}%
}^{+}}\frac{1}{(1-x^{\alpha })}
\end{equation*}
Then $\Delta _{\mathcal{I}}\in \mathcal{F}_{n}$ and $\nabla _{\mathcal{I}%
}\in \mathcal{L}_{n}$. Note that $\nabla _{\mathcal{I}}=\nabla \times \Delta
_{\mathcal{I}}$.

\noindent As a formal series, $\nabla _{\mathcal{I}}$ can be expanded on the
form 
\begin{equation}
\nabla _{\mathcal{I}}=\sum_{\gamma \in \mathbb{Z}^{n}}\mathcal{P}_{\mathcal{I%
}}(\gamma )x^{\gamma }.  \label{expnabI}
\end{equation}
Consider $\lambda \in \mathcal{P}_{n}$ and $\mu =(\mu ^{(0)},...,\mu ^{(p)})$
a dominant weight associated with $G_{\mathcal{I}}$.\ We denote by $%
[V(\lambda ):V_{\mathcal{I}}(\mu )]$ the multiplicity of the irreducible
representation $V_{\mathcal{I}}(\mu )$ in the restriction of $V(\lambda )$
from $G$ to $G_{\mathcal{I}}$.\ The proposition below follows from Theorem
8.2.1 in \cite{GW}:

\begin{proposition}
\label{prop_GW}Consider $\lambda \in \mathcal{P}_{n}$ and $\mu =(\mu
^{(0)},...,\mu ^{(p)})$ a dominant weight of $P_{\mathcal{I}}^{+}$.\ Then 
\begin{equation*}
\lbrack V(\lambda ):V_{\mathcal{I}}(\mu )]=\sum_{w\in W}\varepsilon (w)%
\mathcal{P}_{\mathcal{I}}(w\circ \lambda -\mu ).
\end{equation*}
\end{proposition}

\noindent Define the linear map 
\begin{equation*}
\left\{ 
\begin{tabular}{c}
$\mathrm{S}_{\mathcal{I}}:\mathcal{L}_{n}\rightarrow\widehat{A}$ \\ 
$x^{\beta}\mapsto\mathrm{H}(\Delta_{\mathcal{I}}\times x^{\beta})$%
\end{tabular}
\right. \text{.}
\end{equation*}
For any dominant weight $\mu\in P_{\mathcal{I}}^{+}$, set 
\begin{equation}
S_{\mu,\mathcal{I}}=\mathrm{H}(\Delta_{\mathcal{I}}\times x^{\mu})=\mathrm{S}%
_{\mathcal{I}}(x^{\mu}).  \label{def_Smu}
\end{equation}

\begin{proposition}
\label{prop_dec_Smu}With the above notations we have 
\begin{equation}
S_{\mu ,\mathcal{I}}=\sum_{\lambda \in \mathcal{P}_{n}}[V(\lambda ):V_{%
\mathcal{I}}(\mu )]s_{\lambda }.  \label{dec_S}
\end{equation}
\end{proposition}

\begin{proof}
For any $\beta \in \mathbb{Z}^{n}$, we have obtained in the proof of
Proposition \ref{propHS}, the identity $\mathrm{H}(x^{\beta })=\mathrm{S}%
\circ \chi _{\nabla }(x^{\beta }).$ Denote by $\chi _{\Delta ,\mathcal{I}}$
the linear map defined on $\mathcal{L}_{n}$ by setting $\chi _{\Delta ,%
\mathcal{I}}(x^{\beta })=\Delta _{\mathcal{I}}\times x^{\beta }$. We obtain $%
\mathrm{S}_{\mathcal{I}}(x^{\beta })=\mathrm{H}(\Delta _{\mathcal{I}}\times
x^{\beta })=\mathrm{S}\circ \chi _{\nabla }\circ \chi _{\Delta ,\mathcal{I}%
}(x^{\beta })=\mathrm{S}(\nabla _{\mathcal{I}}\times x^{\beta })$since $%
\nabla _{\mathcal{I}}=\nabla \times \Delta _{\mathcal{I}}$.\ Thus by (\ref
{expnabI}) this yields $\mathrm{S}_{\mathcal{I}}(x^{\beta })=\sum_{\gamma
\in \mathbb{Z}^{n}}\mathcal{P}_{\mathcal{I}}(\gamma )s_{\beta +\gamma }.$
For any $\gamma $, we know by (\ref{SL}) that $s_{\beta +\gamma }=0$ or
there exists $\lambda \in \mathcal{P}_{n},$ $w\in W$ such that $\lambda
=w\circ (\beta +\gamma )$ and $s_{\beta +\gamma }=\varepsilon (w)s_{\lambda
} $.\ This permits to write 
\begin{equation*}
\mathrm{S}_{\mathcal{I}}(x^{\beta })=\sum_{\lambda \in \mathcal{P}%
_{n}}\sum_{w\in W}\varepsilon (w)\mathcal{P}_{\mathcal{I}}(w\circ \lambda
-\beta )s_{\lambda }.
\end{equation*}
When $\beta =\mu $ is a dominant weight of $P_{\mathcal{I}}^{+},$ we obtain
the desired identity by using Proposition \ref{prop_GW}.
\end{proof}

\bigskip

\noindent\textbf{Remarks: }

\noindent $\mathrm{(i):}$ When $G=G_{\mathcal{I}}$ that is, when $r_{0}=n$
and $r_{1}=\cdot \cdot \cdot =r_{p}=0,$ we have $\mu =\mu ^{(0)},$ $\Delta _{%
\mathcal{I}}=\Delta $ and $\mathrm{H}_{\mathcal{I}}=\mathrm{H}$.\ Thus $%
S_{\mu ,\mathcal{I}}=s_{\mu ^{(0)}}$.\ This can be recovered by using (\ref
{dec_S}) since in this case $[V(\lambda ):V_{\mathcal{I}}(\mu )]=0$ except
when $\lambda =\mu ^{(0)}$.

\noindent$\mathrm{(ii):}$ When $G_{\mathcal{I}}=H$ the maximal torus of $G$,
that is when $n=p+1$ and $r_{k}=1$ for any $k=0,...,p,$ we have $%
\mu_{i}=\mu^{(i-1)}$ for any $i=1,...,n$, $\Delta_{\mathcal{I}}=1$ and $%
\mathrm{H}_{\mathcal{I}}(x^{\beta})=h_{\beta}$ for any $\beta\in\mathbb{Z}%
^{n}$. Hence $S_{\mu,\mathcal{I}}=h_{\mu}.\;$In this case $[V(\lambda):V_{%
\mathcal{I}}(\mu)]=K_{\lambda,\mu}$ for any $\lambda\in\mathcal{P}_{n}.$
Thus (\ref{dec_S}) reduces to (\ref{dec_h}).

\noindent$\mathrm{(iii):}$ By $\mathrm{(i)}$ and $\mathrm{(ii)}$ the
functions $S_{\mu,\mathcal{I}}$ interpolates between the Weyl characters $%
s_{\mu}$ and the homogeneous functions $h_{\mu}$.

\noindent $\mathrm{(iv):}$ When $G=GL_{n}$, we have the duality 
\begin{equation*}
\lbrack V(\lambda ):V_{\mathcal{I}}(\mu )]=c_{\mu ^{(0)},...,\mu
^{(p)}}^{\lambda }
\end{equation*}
where $c_{\mu ^{(0)},...,\mu ^{(p)}}^{\lambda }$ is the
Littlewood-Richardson coefficient associated with the multiplicity of $%
V(\lambda )$ in the tensor product $V_{\mu }=V(\mu ^{(0)})\otimes \cdot
\cdot \cdot \otimes V(\mu ^{(p)})$.\ Thus we can write $S_{\mu ,\mathcal{I}%
}=\sum_{\lambda \in \mathcal{P}_{n}}c_{\mu ^{(0)},...,\mu ^{(p)}}^{\lambda
}s_{\lambda }.$ This means that $S_{\mu ,\mathcal{I}}$ is the character of $%
V_{\mu }$. Such a duality does not exist for $G=Sp_{2n},SO_{2n+1}$ or $%
SO_{2n},$ (but see Section \ref{Lastsec}).

\section{Plethysm on Weyl characters\label{Sec-Plethy}}

\subsection{The operators $\Psi_{\ell}$ and $\protect\varphi_{\ell}$}

Consider $\ell$ a positive integer.\ The power sum plethysm operator $%
\Psi_{\ell}$ is defined on $A$ be setting $\Psi_{\ell}(m_{\beta})=m_{\ell%
\beta}$ for any $\beta=(\beta_{1},...,\beta_{n})\in\mathbb{Z}^{n}$ where $%
\ell\beta=(\ell\beta_{1},...,\ell\beta_{n})$. Since $\{m_{\lambda}\mid%
\lambda\in\mathcal{P}_{n}\}$ and $\{h_{\lambda}\mid\lambda\in \mathcal{P}%
_{n}\}$ are dual bases for the inner product $<\cdot,\cdot>$, the adjoint
operator $\varphi_{\ell}$ of $\Psi_{\ell}$ verifies 
\begin{equation}
\varphi_{\ell}(h_{\beta})=\left\{ 
\begin{array}{l}
h_{\beta/\ell}\text{ if }\beta\in(\ell\mathbb{Z})^{n} \\ 
0\text{ otherwise}
\end{array}
\right.  \label{fil_h}
\end{equation}
where $\beta/\ell=(\beta_{1}/\ell,...,\beta_{n}/\ell)$ when $\beta\in (\ell%
\mathbb{Z})^{n}$.

\noindent By abuse of notation, we also denote by $\Psi_{\ell}$ and $%
\varphi_{\ell}$ the linear operators respectively defined on $\mathcal{L}%
_{n} $ by setting 
\begin{equation}
\Psi_{\ell}(x^{\beta})=x^{\ell\beta}\text{ and }\varphi_{\ell}(x^{\beta
})=\left\{ 
\begin{array}{l}
x^{\beta/\ell}\text{ if }\beta\in(\ell\mathbb{Z})^{n} \\ 
0\text{ otherwise}
\end{array}
\right. \text{ for any }\beta\in\mathbb{Z}^{n}\text{.}  \label{fil_ebeta}
\end{equation}

\noindent\textbf{Remark: }Since $\Psi_{\ell}(x^{\beta}\times
x^{\beta^{\prime
}})=\Psi_{\ell}(x^{\beta})\times\Psi_{\ell}(x^{\beta^{\prime}})$ for any $%
\beta,\beta^{\prime}\in\mathbb{Z}^{n}$ the map $\Psi_{\ell}$ is a morphism
of algebra.\ This is not true for $\varphi_{\ell}.\;$Nevertheless, if $%
\{i_{1},...,i_{r}\}$ and $\{j_{1},...,j_{s}\}$ are disjoint subsets of $%
I_{n},$ $\iota=(\iota_{1},...,\iota_{r})\in\mathbb{Z}^{r}$ and $\gamma
=(\gamma_{1},...,\gamma_{s})\in\mathbb{Z}^{s}$ we have 
\begin{equation}
\varphi_{\ell}(x_{i_{1}}^{\iota_{1}}\cdot\cdot\cdot
x_{i_{r}}^{\iota_{r}}\times x_{j_{1}}^{\gamma_{1}}\cdot\cdot\cdot
x_{j_{r}}^{\gamma_{r}})=\varphi_{\ell}(x_{i_{1}}^{\iota_{1}}\cdot\cdot\cdot
x_{i_{r}}^{\iota_{r}})\times\varphi_{\ell}(x_{j_{1}}^{\gamma_{1}}\cdot\cdot%
\cdot x_{j_{r}}^{\gamma_{r}}).  \label{fiprod}
\end{equation}

\bigskip

\noindent For any $\lambda \in \mathcal{P}_{n}$, $\Psi _{\ell }(s_{\lambda
}) $ belongs to $A$, thus decomposes on the basis $\{s_{\mu }\mid \mu \in 
\mathcal{P}_{n}\}.$ Let us write 
\begin{equation*}
\Psi _{\ell }(s_{\lambda })=\sum_{\mu \in \mathcal{P}_{n}}n_{\lambda ,\mu
}s_{\mu }.
\end{equation*}
Since $\Psi _{\ell }$ and $\varphi _{\ell }$ are dual operators with respect
to the scalar product $<\cdot ,\cdot >$, we can write\linebreak\ $n_{\lambda
,\mu }=<\Psi _{\ell }(s_{\lambda }),s_{\mu }>=<s_{\lambda },\varphi _{\ell
}(s_{\mu })>.$ So we have 
\begin{equation*}
\varphi _{\ell }(s_{\mu })=\sum_{\lambda \in \mathcal{P}_{n}}n_{\lambda ,\mu
}s_{\lambda }.
\end{equation*}
By (\ref{S_H}) and Proposition \ref{prop_JT}, we obtain the identity 
\begin{equation*}
s_{\mu }=\sum_{w\in W}\varepsilon (w)h_{\mu +\rho -w\cdot \rho }=\mathrm{H}%
(\Delta \times x^{\mu }).
\end{equation*}
Thus from (\ref{fil_h}) and (\ref{fil_ebeta}) we derive $\varphi _{\ell
}(s_{\mu })=\mathrm{H}(\varphi _{\ell }(\Delta \times x^{\mu })).$ Set $%
P_{\mu }=\Delta \times x^{\mu }.\;$From the previous arguments the
coefficients $n_{\lambda ,\mu }$ are determined by the computation of $%
\varphi _{\ell }(P_{\mu })$.

\subsection{Computation of $\protect\varphi_{\ell}(P_{\protect\mu})\label%
{subsec_comp}$}

For any $i\in \{\overline{n},...,\overline{1}\}$ we set $x_{i}=\frac{1}{x_{%
\overline{i}}}$. This permits to consider also variables indexed by negative
integers. Given $X=(i_{1},...,i_{r})$ an increasing sequence contained in $%
J_{n}$ and $\beta =(\beta _{1},...,\beta _{r})\in \mathbb{Z}^{r}$, we set $%
x_{X}^{\beta }=x_{i_{1}}^{\beta _{1}}\cdot \cdot \cdot x_{i_{r}}^{\beta
_{r}}.$ We also denote by $S_{X}$ the group of permutations of the set $X$.
Each $\sigma \in S_{X}$ determines a unique permutation $\sigma ^{\ast }$ of
the set $\{1,...,r\}$ defined by 
\begin{equation}
\sigma (i_{p})=i_{\sigma ^{\ast }(p)}\text{ for any }p=1,...r.  \label{ide}
\end{equation}
In the sequel, we identify for short $\sigma $ and $\sigma ^{\ast }$.
Similarly, given $Z=(\overline{u}_{r},...,\overline{u}_{1},u_{1},...,u_{r})$
an increasing sequence such that $\{u_{1},...,u_{r}\}\subset I_{n}$, each
signed permutation $w$ defined on $Z$ will be identified with the signed
permutation $w^{\ast }$ defined on $J_{r}$ by $w(u_{p})=u_{w^{\ast }(p)}$
for any $p=1,...r.$

\noindent Set $\rho _{n}=(1,2,...,n)$.\ For any $w\in W$ we have $w\cdot
\rho _{n}=(w(1),...,w(n))$.\ This permits to write 
\begin{equation}
\sum_{w\in W}\varepsilon (w)x^{-w\cdot \rho _{n}}=\sum_{w\in W}\varepsilon
(w)x_{1}^{-w(1)}\cdot \cdot \cdot x_{n}^{-w(n)}.  \label{LRchange}
\end{equation}

\subsubsection{For $G=GL_{n}$}

Set $\kappa _{n}=(1,...,1)\in \mathbb{Z}^{n}$.\ Since $\sigma (\kappa
_{n})=\kappa _{n}$ for any $\sigma \in S_{n}$, one can replace $\rho $ by $%
\rho _{n}=(1,2,...,n)$ in (\ref{iden_delta}). By using (\ref{LRchange}) we
can write 
\begin{equation*}
P_{\mu }=x_{1}^{(\mu _{1}+1)}\cdot \cdot \cdot x_{n}^{(\mu
_{n}+n)}\sum_{\sigma \in S_{n}}\varepsilon (\sigma )x_{1}^{-\sigma (1)}\cdot
\cdot \cdot x_{n}^{-\sigma (n)}
\end{equation*}
where $\mu =(\mu _{1},...,\mu _{n}).\;$For any $k\in \{0,...,\ell -1\}$
consider the ordering sequences 
\begin{equation}
I^{(k)}=(i\in I_{n}\mid \mu _{i}+i\equiv k\mathrm{mod}\ell )\text{ and }%
J^{(k)}=(i\in I_{n}\mid i\equiv k\mathrm{mod}\ell ).  \label{defIJ}
\end{equation}
Set $r_{k}=\mathrm{card}(I^{(k)})$ and write $%
I^{(k)}=(i_{1}^{(k)},...,i_{r_{k}}^{(k)})$. Then 
\begin{equation*}
\mu ^{(k)}=\left( \frac{\mu _{i}+i+\ell -k}{\ell }\mid i\in I^{(k)}\right)
\in \mathbb{Z}^{r_{k}}
\end{equation*}
We derive 
\begin{equation*}
P_{\mu }=x_{I^{(0)}}^{\ell \mu ^{(0)}}x_{I^{(1)}}^{\ell \mu ^{(1)}}\cdot
\cdot \cdot x_{I^{(\ell -1)}}^{\ell \mu ^{(\ell -1)}}\sum_{\sigma \in
S_{n}}\varepsilon (\sigma )x_{1}^{-\sigma (1)}\cdot \cdot \cdot
x_{n}^{-\sigma (n)}\times \prod_{k=0}^{\ell
-1}\prod_{a=1}^{r_{k}}x_{i_{a}^{(k)}}^{-(\ell -k)}.
\end{equation*}
This gives 
\begin{equation}
\varphi _{\ell }(P_{\mu })=x_{I^{(0)}}^{\mu ^{(0)}}x_{I^{(1)}}^{\mu
^{(1)}}\cdot \cdot \cdot x_{I^{(\ell -1)}}^{\mu ^{(\ell -1)}}\sum_{\sigma
\in S_{n}}\varepsilon (\sigma )\varphi _{\ell }\left( \prod_{k=0}^{\ell
-1}\prod_{a=1}^{r_{k}}x_{i_{a}^{(k)}}^{-\sigma (i_{a}^{(k)})-(\ell
-k)}\right) .  \label{eq1}
\end{equation}
The contribution of a fixed permutation $\sigma \in S_{n}$ in the above sum
is nonzero if and only if for any $k=0,...,\ell -1$%
\begin{equation*}
i\in I^{(k)}\Longrightarrow \sigma (i)\equiv k\mathrm{mod}\ell .
\end{equation*}
Thus we must have $\sigma (I^{(k)})\subset J^{(k)}$ for any $k=0,...,\ell
-1. $ Since $\sigma $ is a bijection, $I^{(k)}\cap I^{(k^{\prime
})}=J^{(k)}\cap J^{(k^{\prime })}=\emptyset $ if $k\neq k^{\prime }$and $%
\cup _{0\leq k\leq \ell -1}I^{(k)}=\cup _{0\leq k\leq \ell -1}J^{(k)}=I_{n}$%
, the restriction of $\sigma $ on $I_{k}$ is a bijection from $I^{(k)}$ to $%
J^{(k)}.$ In particular $\mathrm{card}(J^{(k)})=\mathrm{card}%
(I^{(k)})=r_{k}. $ This means that we have the equivalences 
\begin{equation}
\varphi _{\ell }\left( \prod_{k=0}^{\ell
-1}\prod_{a=1}^{r_{k}}x_{i_{a}^{(k)}}^{-\sigma (i_{a}^{(k)})-(\ell
-k)}\right) \neq 0\Longleftrightarrow \sigma (I^{(k)})=J^{(k)}\text{ for any 
}k=0,...,\ell -1.  \label{cond}
\end{equation}
Write $J^{(k)}=(k,k+\ell ,...,k+(r_{k}-1)\ell ).$ Denote by $\sigma _{0}\in
S_{n}$ the permutation verifying 
\begin{equation}
\sigma _{0}(i_{a}^{(k)})=k+(a-1)\ell  \label{def_sigma0}
\end{equation}
for any $k=0,...,\ell -1$ and any $a=1,...,r_{k}.$ Let $S_{I^{(k)}}$ be the
permutation group of the set $I^{(k)}.$ The permutations $\sigma $ which
verify the right-hand side of (\ref{cond}) can be written $\sigma =\sigma
_{0}\tau $ where $\tau =(\tau ^{(0)},\cdot \cdot \cdot ,\tau ^{(\ell -1)})$
belongs to the direct product $S_{I^{(0)}}\times \cdot \cdot \cdot \times
S_{I^{(\ell -1)}}$.\ We have then $\varepsilon (\sigma )=\varepsilon (\sigma
_{0})(-1)^{l(\tau ^{(0)})}\times \cdot \cdot \cdot \times (-1)^{l(\tau
^{(p)})}$

\noindent For any $k\in \{0,...,\ell -1\}$ set 
\begin{equation*}
P_{k}=\sum_{\tau ^{(k)}\in S_{I^{(k)}}}(-1)^{l(\tau )}\varphi _{\ell }\left(
\prod_{a=1}^{r_{k}}x_{i_{a}^{(k)}}^{-\sigma _{0}\tau
^{(k)}(i_{a}^{(k)})-(\ell -k)}\right) .
\end{equation*}
From (\ref{fiprod}) and (\ref{eq1}) we derive 
\begin{equation*}
\varphi _{\ell }(P_{\mu })=\varepsilon (\sigma _{0})\prod_{k=0}^{\ell
-1}x_{I^{(k)}}^{\mu ^{(k)}}P_{k}.
\end{equation*}
Since $\sigma _{0}(i_{a}^{(k)})=k+(a-1)\ell $, we can write by (\ref{ide}) $%
\sigma _{0}\tau ^{(k)}(i_{a}^{(k)})=k+(\tau ^{(k)}(a)-1)\ell $. Thus we
obtain 
\begin{equation*}
P_{k}=\sum_{\tau ^{(k)}\in S_{I^{(k)}}}(-1)^{l(\tau )}x_{i_{1}^{(k)}}^{-\tau
^{(k)}(1)}\cdot \cdot \cdot x_{i_{r_{k}}^{(k)}}^{-\tau
^{(k)}(r_{k})}=x_{I^{(k)}}^{-\rho _{r_{k}}}\Delta _{I^{(k)}}
\end{equation*}
where $\rho _{r_{k}}=(1,2,...,r_{k})$ and $\Delta _{I^{(k)}}=\prod_{i<j\
i,j\in I^{(k)}}(1-x_{j}/x_{i}).$ Finally, this gives 
\begin{equation*}
\varphi _{\ell }(P_{\mu })=\varepsilon (\sigma _{0})\prod_{k=0}^{\ell
-1}x_{I^{(k)}}^{\mu ^{(k)}-\rho _{r_{k}}}\Delta _{I^{(k)}}=\varepsilon
(\sigma _{0})\prod_{k=0}^{\ell -1}x_{I^{(k)}}^{\mu ^{(k)}}\Delta _{I^{(k)}}
\end{equation*}
where for any $k=0,...,\ell -1,$%
\begin{equation}
\mu ^{(k)}=\left( \frac{\mu _{i}+i+\ell -k}{\ell }\mid i\in I^{(k)}\right)
-(1,2,...,r_{k})\in \mathbb{Z}^{r_{k}}.  \label{def_muk}
\end{equation}

\begin{theorem}
\label{thA}Consider a partition $\mu $ of length $n$ and $\ell $ a positive
integer.$\;$For any $k=0,...,\ell -1$ define the sets $I^{(k)}$ and $J^{(k)}$
as in (\ref{defIJ}).

\begin{itemize}
\item  If there exists $k\in \{0,...,\ell -1\}$ such that $\mathrm{card}%
(I^{(k)})\neq \mathrm{card}(J^{(k)})$ then $\varphi _{\ell }(s_{\mu })=0.$

\item  Otherwise, for any $k=0,...,\ell -1,$ set $r_{k}=\mathrm{card}%
(I^{(k)})=\mathrm{card}(J^{(k)})$ and define $\sigma _{0}$ as in (\ref
{def_sigma0}).\ Then each $r_{k}$-tuple defined by (\ref{def_muk}) is a
partition and we have 
\begin{equation*}
\varphi _{\ell }(s_{\mu })=\varepsilon (\sigma _{0})S_{\binom{\mu }{\ell },%
\mathcal{I}}=\varepsilon (\sigma _{0})\mathrm{char}(V_{\mu })
\end{equation*}
where $\mathcal{I=}\{I^{(0},...,I^{(\ell -1)}\},$ $\binom{\mu }{\ell }=(\mu
^{(0)},...,\mu ^{(\ell -1)})$ and $\mathrm{char}(V_{\mu })$ is the character
of the $GL_{n}$-module $V_{\mu }=V(\mu ^{(0)})\otimes \cdot \cdot \cdot
\otimes V(\mu ^{(\ell -1)}).$
\end{itemize}
\end{theorem}

\begin{proof}
One verifies easily from (\ref{def_muk}) that each $\mu ^{(k)}$ is a
partition. By the previous computation, we obtain 
\begin{equation*}
\varphi _{\ell }(\Delta \times x^{\mu })=\varepsilon (\sigma _{0})\Delta _{%
\mathcal{I}}\times x^{\binom{\mu }{\ell }}
\end{equation*}
(with the notation of \ref{subsec_BC}). By definition of $\varphi _{\ell }$
we have also 
\begin{equation*}
\varphi _{\ell }(s_{\mu })=\varepsilon (\sigma _{0})\mathrm{H}\circ \varphi
_{\ell }(\Delta \times x^{\mu })=\varepsilon (\sigma _{0})\mathrm{H}(\Delta
_{\mathcal{I}}\times x^{\binom{\mu }{\ell }})=\varepsilon (\sigma _{0})S_{%
\binom{\mu }{\ell },\mathcal{I}}
\end{equation*}
where the last equality follows from (\ref{def_Smu}).
\end{proof}

\bigskip

\noindent \textbf{Remark: }The subgroup $G_{\mathcal{I}}$ appearing in
Theorem \ref{thA} is characterized by $\mathcal{I=}\{I^{(0},...,I^{(\ell
-1)}\}.$ This means that for type $A$, we have $X^{(k)}=I^{(k)}$ for any $%
k>0 $ with the notation of \ref{subsec}, that is the sets $X^{(k)}$ contain
only positive indices.

\begin{example}
\ \ 

\noindent Consider $\mu =(1,2,3,4,4,4,6,6)$ and take $\ell =3.$ We have $\mu
+\rho _{8}=(2,4,6,8,9,10,13,14).$ Thus $I^{(0)}=\{3,5\},I^{(1)}=\{2,6,7%
\},I^{(2)}=\{1,4,8\}$ and $J^{(0)}=\{3,6\},J^{(1)}=\{1,4,7\},J^{(2)}=\{2,5,8%
\}.$ Then $\mu ^{(0)}=(1,1),\mu ^{(1)}=(1,2,2)$ and $\mu ^{(2)}=(0,1,2).$ We
have $G_{\mathcal{I}}\simeq GL_{2}\times GL_{3}\times GL_{3}.$
\end{example}

\subsubsection{For $G=Sp_{2n}\label{subsubSP}$}

We have $\rho =\rho _{n}=(1,2,...,n)$. By using (\ref{LRchange}) we deduce
the identity: 
\begin{equation}
P_{\mu }=x_{1}^{(\mu _{1}+1)}\cdot \cdot \cdot x_{n}^{(\mu
_{n}+n)}\sum_{w\in W}\varepsilon (w)x_{1}^{-w(1)}\cdot \cdot \cdot
x_{n}^{-w(n)}  \label{pmuC}
\end{equation}
where $W$ is the group of signed permutations defined on $J_{n}=\{\overline{n%
},...,\overline{1},1,...,n\},$ that is the subgroup of permutations $w\in
S_{J_{n}}$ verifying $w(\overline{x})=\overline{w(x)}$ for any $x\in J_{n}$%
.\ Given $k\in \{0,...,\ell -1\}$ consider the ordering sequences 
\begin{equation}
I^{(k)}=(i\in I_{n}\mid \mu _{i}+i\equiv k\mathrm{mod}\ell )\text{ and }%
J^{(k)}=(x\in J_{n}\mid x\equiv k\mathrm{mod}\ell ).  \label{deIJC}
\end{equation}
Set $p=\frac{\ell }{2}$ if $\ell $ is even and $p=\frac{_{\ell -1}}{2}$
otherwise.

\paragraph{The odd case $\ell=2p-1$}

Set $r_{0}=\mathrm{card}(I^{(0)})$ and for any $k=1,...,p-1$, $s_{k}=\mathrm{%
card}(I_{k}),$ $r_{k}=\mathrm{card}(I_{k})+\mathrm{card}(I_{\ell-k})$.\
Write $X^{(k)},k=1,...,p$ for the increasing reordering of $\overline
{I}_{k}\cup I_{\ell-k}.$ Set $I^{(0)}=(i_{1}^{(0)},...,i_{r_{0}}^{(0)})$ and
for $k>0$%
\begin{equation}
X^{(k)}=(i_{1}^{(k)},...,i_{r_{k}}^{(k)}).  \label{XKC}
\end{equation}
This means that $I^{(k)}=(\overline{i}_{s_{k}}^{(k)},...,\overline{i}%
_{1}^{(k)})$ and $I^{(\ell-k)}=(i_{s_{k+1}}^{(k)},...,i_{r_{k}}^{(k)})$. To
simplify the computation, we are going to use the indices and the variables $%
x_{i},i\in X^{(k)}$ rather than the variables $x_{i},i\in I^{(k)}\cup
I^{(\ell-k)}$ when $k\in\{1,...,p-1\}.$

\noindent Consider 
\begin{equation*}
\mu ^{(0)}=\left( \frac{\mu _{i}+i}{\ell }\mid i\in I^{(0)}\right) \in 
\mathbb{Z}^{r_{0}}\text{ and for }k>0,\text{ }\mu ^{(k)}=\left( \mathrm{sign}%
(i)\frac{\mu _{\left| i\right| }+\left| i\right| +\mathrm{sign}(i)k}{\ell }%
\mid i\in X^{(k)}\right) \in \mathbb{Z}^{r_{k}}
\end{equation*}
where for any $i\in J_{n}$, $\mathrm{sign}(i)=1$ if $i>0$ and $-1$
otherwise.\ For any $i\in I$, we have $x_{i}^{-w(i)}=x_{\overline{i}}^{-w(%
\overline{i})}.$ Thus 
\begin{equation*}
\prod_{i\in X^{(k)}}x_{i}^{-w(i)}=\prod_{i\in
I^{(k)}}x_{i}^{-w(i)}\prod_{i\in I^{(\ell -k)}}x_{i}^{-w(i)}\text{ and }%
x_{1}^{(\mu _{1}+1)}\cdot \cdot \cdot x_{n}^{(\mu _{n}+n)}=x_{I^{(0)}}^{\ell
\mu ^{(0)}}\prod_{k=1}^{p-1}x_{X^{(k)}}^{\ell \mu ^{(k)}}\prod_{i\in
X^{(k)}}x_{i}^{-k}
\end{equation*}
by definition of the $\mu ^{(k)}$'s. Then (\ref{pmuC}) can be rewritten\ 
\begin{equation*}
P_{\mu }=x_{I^{(0)}}^{\ell \mu ^{(0)}}\prod_{k=1}^{p-1}x_{X^{(k)}}^{\ell \mu
^{(k)}}\times \sum_{w\in W}\varepsilon (w)\prod_{i\in
I^{(0)}}x_{i}^{-w(i)}\prod_{k=1}^{p-1}\prod_{i\in
X^{(k)}}x_{i}^{-w(i)}\times \prod_{k=1}^{p-1}\prod_{i\in X^{(k)}}x_{i}^{-k}.
\end{equation*}
This gives 
\begin{equation*}
\varphi _{\ell }(P_{\mu })=x_{I^{(0)}}^{\mu
^{(0)}}\prod_{k=1}^{p-1}x_{X^{(k)}}^{\mu ^{(k)}}\times \sum_{w\in
W}\varepsilon (w)\varphi _{\ell }\left( \prod_{i\in
I^{(0)}}x_{i}^{-w(i)}\prod_{k=1}^{p-1}\prod_{i\in
X^{(k)}}x_{i}^{-w(i)-k}\right)
\end{equation*}
The contribution of a fixed $w\in W$ in the above sum is nonzero if and only
if 
\begin{equation}
\left\{ 
\begin{array}{l}
i\in I^{(0)}\Longrightarrow w(i)\equiv 0\mathrm{mod}\ell \text{ } \\ 
i\in X^{(k)}\Longrightarrow w(i)\equiv -k\mathrm{mod}\ell \text{ for any }%
k=1,...,p-1
\end{array}
\right.  \label{condC}
\end{equation}
Thus we must have $w(\overline{I}^{(0)}\cup I^{(0)})\subset J^{(0)}$ and for
any $k=1,...,p-1,$ $w(X^{(k)})\subset J^{(\ell -k)}$.\ Recall that $%
\overline{J}^{(0)}=J^{(0)}$ and $\overline{J}^{(\ell -k)}=J^{(k)}$ for $%
k=1,...,p-1.$ Moreover 
\begin{equation*}
I^{(0)}\cup \overline{I}^{(0)}\bigcup_{k=1}^{p-1}X^{(k)}\cup \overline{X}%
^{(k)}=J_{n}\text{ and }J^{(0)}\cup \bigcup_{k=1}^{p-1}J^{(k)}\cup J^{(\ell
-k)}=J_{n}.
\end{equation*}
Since the sets appearing in the left hand side of these two equalities are
pairwise disjoint, we must have $w(\overline{I}^{(0)}\cup I^{(0)})=J^{(0)},$
and for $k=1,...,p-1,$ $w(X^{(k)})=J^{(\ell -k)}$.\ In particular $\mathrm{%
card}(J^{(0)})=2\mathrm{card}(I^{(0)})=2r_{0}$ and $\mathrm{card}(J^{(\ell
-k)})=\mathrm{card}(X^{(k)})=r_{k}$ for any $k=1,...,p-1.\;$We have the
equivalences 
\begin{equation}
\varphi _{\ell }\left( \prod_{i\in
I^{(0)}}x_{i}^{-w(i)}\prod_{k=1}^{p-1}\prod_{i\in
X^{(k)}}x_{i}^{-w(i)+k}\right) \neq 0\Longleftrightarrow \left\{ 
\begin{tabular}{l}
$\mathrm{(i):}$ $w(I^{(0)}\cup \overline{I}^{(0)})=J^{(0)}$ \\ 
$\mathrm{(ii):}$ $w(X^{(k)})=J^{(\ell -k)}$ for any $k=1,...,p-1$%
\end{tabular}
\right. .  \label{conC}
\end{equation}
Note that condition $\mathrm{(ii)}$ can be rewritten: $w(\overline{X}%
^{(k)})=J^{(k)}$ for any $k=1,...,p-1.\;$

\noindent We can set $J^{(0)}=(-r_{0}\ell,...,r_{0}\ell)$ and for $%
k=1,...,p-1,$ 
\begin{equation*}
J^{(\ell-k)}=(-k-\alpha_{k}\ell,...,-k+\beta_{k}\ell)\text{, }%
J^{(k)}=(k-\beta_{k}\ell,...,k+\alpha_{k}\ell)
\end{equation*}
with $\alpha_{k}+\beta_{k}+1=r_{k}.$

\noindent Consider $w_{0}\in W$ defined by 
\begin{align}
w_{0}(i_{a}^{(0)})& =a\ell \text{ for }a\in \{1,...,r_{0}\}  \label{w0C} \\
w_{0}(i_{a}^{(k)})& =-k-\alpha _{k}\ell +(a-1)\ell \text{ for any }%
k=1,...,p-1.  \notag
\end{align}
Denote by $\mathcal{W}$ the set of signed permutations $w$ which verify $%
\mathrm{(i)}$ and $\mathrm{(ii)}$ in (\ref{conC}). We have $w_{0}\in 
\mathcal{W}$.\ Each $w\in \mathcal{W}$ can be written $w=w_{0}v$ where $%
v=(v^{(0)},\tau ^{(1)},...,\tau ^{(p-1)})$ belongs to the direct product $%
W_{I^{(0)}}\times S_{X^{(1)}}\times \cdot \cdot \cdot \times S_{X^{(p-1)}}$.
Here $W_{I^{(0)}}$ is the group of signed permutations defined on $\overline{%
I}^{(0)}\cup I^{(0)}$ and for $k=1,...,p-1,$ $S_{X^{(k)}}$ is the group of
signed permutations $\tau ^{(k)}$ defined on $\overline{X}^{(k)}\cup X^{(k)}$
and verifying $\tau ^{(k)}(X^{(k)})=X^{(k)}$. Indeed if $\tau ^{(k)}(x)\in 
\overline{X}^{(k)}$ and $x\in X^{(k)}$, we would have $w(x)\in J^{(k)}$ and $%
x\in X^{(k)}$ which contradicts $\mathrm{(ii)}$. This means that $%
S_{X^{(k)}} $ is in fact isomorphic to the symmetric group $S_{r_{k}}.$
Since the sets $I^{(0)}$ and $X^{(k)},k=1,...,p-1$ are increasing
subsequences of $J_{n}$, we have by Lemma \ref{lem_sign} $\varepsilon
(w)=\varepsilon (w_{0})(-1)^{l(v^{(0)})}(-1)^{l(\tau ^{(1)})}\times \cdot
\cdot \cdot \times (-1)^{l(\tau ^{(p-1)})}$.

\noindent Set 
\begin{align*}
P_{0}& =\sum_{v^{(0)}\in W_{I^{(0)}}}(-1)^{l(v^{(0)})}\varphi _{\ell }\left(
\prod_{i\in I^{(0)}}x_{i}^{-w_{0}v^{(0)}(i)}\right) \text{ and} \\
P_{k}& =\sum_{\tau ^{(k)}\in S_{X^{(k)}}}(-1)^{l(\tau ^{(k)})}\varphi _{\ell
}\left( \prod_{i\in X^{(k)}}x_{i}^{-w_{0}\tau ^{(k)}(i)-k}\right) ,k\in
\{1,...,p-1\}.
\end{align*}
We obtain 
\begin{equation*}
\varphi _{\ell }(P_{\mu })=\varepsilon (w_{0})x_{I^{(0)}}^{\mu
^{(0)}}P_{0}\prod_{k=1}^{p-1}x_{X^{(k)}}^{\mu ^{(k)}}P_{k}.
\end{equation*}
From (\ref{ide}) and (\ref{w0C}), we have $%
w_{0}v^{(0)}(i_{a}^{(0)})=v^{(0)}(a)\ell $ for any $a=1,...,r_{0}$ and 
\begin{equation*}
w_{0}\tau ^{(k)}(i_{a}^{(k)})=-k-\alpha _{k}\ell +(\tau ^{(k)}(a)-1)\ell 
\text{ for any }a=1,...,r_{k}.
\end{equation*}
This yields 
\begin{align*}
P_{0}& =\sum_{v^{(0)}\in
W_{I^{(0)}}}(-1)^{l(v^{(0)})}%
\prod_{a=1}^{r_{0}}x_{i_{a}^{(0)}}^{-v^{(0)}(a)}=x_{I^{(0)}}^{-\rho
_{r_{0}}}\Delta _{I^{(0)}}\text{ and } \\
P_{k}& =\sum_{\tau ^{(k)}\in S_{X^{(k)}}}(-1)^{l(\tau
^{(k)})}\prod_{a=1}^{r_{k}}x_{i_{a}^{(k)}}^{-\tau ^{(k)}(a)+(\alpha
_{k}+1)}=x_{X^{(k)}}^{\eta _{r_{k}}}\Delta _{X^{(k)}}
\end{align*}
where for any $k=1,...,p-1,$ $\eta _{r_{k}}=-\rho _{r_{k}}+(\alpha
_{k}+1,...,\alpha _{k}+1)\in \mathbb{Z}^{r_{k}},$%
\begin{align*}
\Delta _{I^{(0)}}& =\prod_{i<j\ i,j\in I^{(0)}}(1-\frac{x_{j}}{x_{i}}%
)\prod_{r\leq s\ r,s\in I^{(0)}}(1-x_{r}x_{s}) \\
\text{and }\Delta _{X^{(k)}}& =\prod_{i<j\ i,j\in X^{(k)}}(1-\frac{x_{j}}{%
x_{i}})\text{ for any }k=1,...,p-1.
\end{align*}
Finally, this gives 
\begin{equation*}
\varphi _{\ell }(P_{\mu })=\varepsilon (w_{0})x_{I^{(0)}}^{\mu ^{(0)}-\rho
_{r_{0}}}\Delta _{I^{(0)}}\prod_{k=1}^{p-1}x_{X^{(k)}}^{\mu ^{(k)}-\eta
_{r_{k}}}\Delta _{X^{(k)}}=\varepsilon (w_{0})x_{I^{(0)}}^{\mu ^{(0)}}\Delta
_{I^{(0)}}\prod_{k=1}^{p-1}x_{X^{(k)}}^{\mu ^{(k)}}\Delta _{X^{(k)}}
\end{equation*}
where 
\begin{equation}
\mu ^{(0)}=\left( \frac{\mu _{i}+i}{\ell }\mid i\in I^{(0)}\right)
-(1,...,r_{0})\in \mathbb{Z}^{r_{0}}  \label{mu0C}
\end{equation}
and for any $k=1,...,p-1$%
\begin{equation}
\mu ^{(k)}=\left( \mathrm{sign}(i)\frac{\mu _{\left| i\right| }+\left|
i\right| +\mathrm{sign}(i)k}{\ell }\mid i\in X^{(k)}\right)
-(1,...,r_{k})+(\alpha _{k}+1,....,\alpha _{k}+1)\in \mathbb{Z}^{r_{k}}.
\label{mukC}
\end{equation}
Recall that the weights corresponding to the subgroup of Levi type $G_{%
\mathcal{I}}$ are written following the convention (\ref{conv}).

\begin{theorem}
\label{thC}Consider a partition $\mu $ of length $n$ and $\ell =2p-1$ a
positive integer.$\;$Let $I^{(0)}$ and $J^{(0)}$ be as in (\ref{deIJC}).\
For any $k=1,...,p-1$ define the sets $X^{(k)}$ and $J^{(k)}$ by (\ref{deIJC}%
) and (\ref{XKC}).

\begin{itemize}
\item  If $\mathrm{card}(I^{(0)})\neq \frac{1}{2}\mathrm{card}(J^{(0)})$ or
if there exists $k\in \{1,...,p-1\}$ such that $\mathrm{card}(X^{(k)})\neq 
\mathrm{card}(J^{(k)})$ then $\varphi _{\ell }(s_{\mu })=0.$

\item  Otherwise, set $r_{0}=\mathrm{card}(I^{(0)})$ and for any $%
k=1,...,p-1,$ $r_{k}=\mathrm{card}(X^{(k)})$.\ Let $w_{0}\in W$ be as in (%
\ref{w0C}).\ Consider $\binom{\mu }{\ell }=(\mu ^{(0)},\mu ^{(1)},...,\mu
^{(p-1)})$ where the $\mu ^{(k)}$'s are defined by (\ref{mu0C}) and (\ref
{mukC}). Then $\binom{\mu }{\ell }$ is a dominant weight of $P_{\mathcal{I}%
}^{+}$ with $\mathcal{I=}\{I^{(0)},X^{(1)}...,X^{(p-1)}\}$ and we have 
\begin{equation*}
\varphi _{\ell }(s_{\mu })=\varepsilon (w_{0})S_{\binom{\mu }{\ell },%
\mathcal{I}}.
\end{equation*}
\end{itemize}
\end{theorem}

\begin{proof}
The proof is essentially the same as in Theorem \ref{thA}. We obtain 
\begin{equation*}
\varphi _{\ell }(\Delta \times x^{\mu })=\varepsilon (w_{0})\Delta _{%
\mathcal{I}}\times x^{\binom{\mu }{\ell }}
\end{equation*}
where in the right-hand side of the preceding equality $\binom{\mu }{\ell }$
is expressed on the basis $\{\varepsilon _{i}\mid i\in I_{n}\}$ (see (\ref
{iden})). This permits to write as in the case $G=GL_{n}$ 
\begin{equation*}
\varphi _{\ell }(s_{\mu })=\varepsilon (w_{0})\mathrm{H}\circ \varphi _{\ell
}(\Delta \times x^{\mu })=\varepsilon (w_{0})\mathrm{H}(\Delta _{\mathcal{I}%
}\times x^{\binom{\mu }{\ell }})=\varepsilon (w_{0})S_{\binom{\mu }{\ell },%
\mathcal{I}}.
\end{equation*}
\end{proof}

\begin{example}
\ \ 

\noindent Consider $\mu =(1,2,3,4,4,4,6,6)$ and take $\ell =3.$ We have $\mu
+\rho _{8}=(2,4,6,8,9,10,13,14).$ Thus $I^{(0)}=\{3,5\},X^{(1)}=\{\overline{7%
},\overline{6},\overline{2},1,4,8\}$ and $J^{(0)}=\{\overline{6},\overline{3}%
,3,6\},J^{(1)}=\{\overline{8},\overline{5},\overline{2},1,4,7\},J^{(2)}=\{%
\overline{7},\overline{4},\overline{1},2,5,8\}.$ In particular $\alpha
_{1}=2.\;$Then $\mu ^{(0)}=(1,1)$ and $\mu ^{(1)}=$%
\begin{multline*}
\left( -\frac{13-1}{3}-1+3,-\frac{10-1}{3}-2+3,-\frac{4-1}{3}-3+3,\frac{2+1}{%
3}-4+3,\frac{8+1}{3}-5+3,\frac{14+1}{3}-6+3\right) \\
=(-2,-2,-1,0,1,2)
\end{multline*}
with the convention (\ref{conv}). We have $G_{\mathcal{I}}\simeq
Sp_{4}\times GL_{6}.$
\end{example}

\paragraph{The even case $\ell=2p$}

With the same notation as in the odd case, (\ref{pmuC}) can be rewritten\ 
\begin{equation*}
P_{\mu }=x_{I^{(0)}}^{\ell \mu ^{(0)}}x_{I^{(p)}}^{\ell \mu
^{(p)}}\prod_{k=1}^{p-1}x_{X^{(k)}}^{\ell \mu ^{(k)}}\times \sum_{w\in
W}\varepsilon (w)\prod_{i\in I^{(0)}}x_{i}^{-w(i)}\prod_{i\in
I^{(p)}}x_{i}^{-w(i)-p}\prod_{k=1}^{p-1}\prod_{i\in X^{(k)}}x_{i}^{-w(i)-k}.
\end{equation*}
where 
\begin{equation*}
\mu ^{(p)}=\left( \frac{\mu _{i}+i+p}{\ell }\mid i\in I^{(p)}\right) .
\end{equation*}
This gives 
\begin{equation*}
\varphi _{\ell }(P_{\mu })=x_{I^{(0)}}^{\mu ^{(0)}}x_{I^{(p)}}^{\mu
^{(p)}}\prod_{k=1}^{p-1}x_{X^{(k)}}^{\mu ^{(p)}}\times \sum_{w\in
W}\varepsilon (w)\varphi _{\ell }\left( \prod_{i\in
I^{(0)}}x_{i}^{-w(i)}\prod_{i\in
I^{(p)}}x_{i}^{-w(i)-p}\prod_{k=1}^{p-1}\prod_{i\in
X^{(k)}}x_{i}^{-w(i)-k}\right)
\end{equation*}
The contribution of a fixed $w\in W$ in the above sum is nonzero if
conditions (\ref{conC}) are verified and 
\begin{equation*}
i\in I^{(p)}\Longrightarrow w(i)\in J^{(p)}.
\end{equation*}
Since $p\equiv -p\mathrm{mod}\ell $ we have $J^{(p)}=\overline{J}%
^{(p)}=\{-p-\alpha _{p}\ell ,...,-p,p,...,p+\alpha _{p}\ell \}$. This
implies that $w(I^{(p)}\cup \overline{I}^{(p)})=J^{(p)}$ and thus $\mathrm{%
card}(I^{(p)})=\frac{1}{2}\mathrm{card}(J^{(p)}).\;$We then define $w_{0}$
by requiring (\ref{w0C}) and $w_{0}(i_{a}^{(p)})=p+(a-1)\ell $ for $a\in
\{1,...,r_{p}\}.\;$By using similar arguments as in the odd case, we obtain
that $w$ can be written $w=w_{0}\nu $ where $\nu =(v^{(0)},\tau
^{(1)},...,\tau ^{(p-1)},v^{(p})$ belongs to the direct product $%
W_{I^{(0)}}\times S_{X^{(1)}}\times \cdot \cdot \cdot \times S_{X^{(\ell
-1)}}\times W_{I^{(p)}}$ with $W_{I^{(p)}}$ the group of signed permutations
defined on $\overline{I^{(p)}}\cup I^{(p)}$.\ Note that $W_{I^{(p)}}$ is a
Weyl group of type $B_{r_{p}}.$ By Lemma \ref{lem_sign} we have also $%
\varepsilon (w)=\varepsilon (w_{0})(-1)^{l(v^{(0)})}(-1)^{l(\tau
^{(1)})}\times \cdot \cdot \cdot \times (-1)^{l(\tau ^{(p)})}\times
(-1)^{l(v^{(p)})}$. We obtain 
\begin{equation*}
\varphi _{\ell }(P_{\mu })=\varepsilon (w_{0})x_{I^{(0)}}^{\mu
^{(0)}}P_{0}\times x_{I^{(p)}}^{\mu
^{(p)}}P_{p}\prod_{k=1}^{p-1}x_{X^{(k)}}^{\mu ^{(k)}}P_{k}\text{ where }%
P_{p}=\sum_{v^{(p)}\in W_{I^{(p)}}}(-1)^{l(v^{(p)})}\varphi _{\ell }\left(
\prod_{i\in I^{(p)}}x_{i}^{-w_{0}v^{(p)}(i)-p}\right) .
\end{equation*}
The functions $P_{k},k=0,...,p-1$ can be computed as in the odd case.\ For $%
P_{p}$, observe that each $v^{(p)}\in W_{I^{(p)}}$ can be written $%
v^{(p)}=\zeta \sigma $ according to the decomposition of $W_{I^{(p)}}$ as
the semidirect product $(\mathbb{Z}/2\mathbb{Z)}^{r_{p}}\varpropto
S_{I^{(p)}}$. We have then for any $a=1,...,r_{p}$, $%
w_{0}v^{(p)}(i_{a}^{(p)})=\xi (a)(p+(\sigma (a)-1)\ell ).$ This yields 
\begin{equation*}
P_{p}=\sum_{v^{(p)}\in W_{I^{(p)}}}(-1)^{l(v^{(p)})}\varphi _{\ell }\left(
\prod_{a=1}^{r_{p}}x_{i}^{-\xi (a)(p+(\sigma (a)-1)\ell )-p}\right)
=\sum_{v^{(p)}\in
W_{I^{(p)}}}(-1)^{l(v^{(p)})}\prod_{a=1}^{r_{p}}x_{i_{a}^{(p)}}^{-\frac{%
1-\xi (a)}{2}-\xi \sigma (a)}.
\end{equation*}
Thus 
\begin{equation*}
P_{p}=\prod_{i\in I^{(p)}}x_{i}^{-1/2}\sum_{v^{(p)}\in
W_{I^{(p)}}}(-1)^{l(v^{(p)})}\left( v^{(p)}\cdot
\prod_{a=1}^{r_{p}}x_{i_{a}^{(p)}}^{-(a-\frac{1}{2})}\right)
=x_{I^{(p)}}^{-\rho _{r_{p}}}\Delta _{I^{(p)},B_{r_{p}}}
\end{equation*}
where 
\begin{equation*}
\Delta _{I^{(p)},B_{r_{p}}}=\prod_{i<j\ i,j\in I^{(p)}}(1-\frac{x_{j}}{x_{i}}%
)\prod_{r<s\ r,s\in I^{(p)}}(1-x_{r}x_{s})\prod_{i\in I^{(p)}}(1-x_{i}).
\end{equation*}
Indeed the half sum of positive roots is equal to $(\frac{1}{2},...,r_{p}-%
\frac{1}{2})$ in type $B_{r_{p}}$.\ This means that when $\ell $ is even 
\begin{equation*}
\varphi _{\ell }(P_{\mu })=\varepsilon (w_{0})x_{I^{(0)}}^{\mu ^{(0)}}\Delta
_{I^{(0)}}\prod_{k=1}^{p-1}x_{X^{(k)}}^{\mu ^{(k)}}\Delta _{X^{(k)}}\times
x_{I^{(p)}}^{\mu ^{(p)}}\Delta _{I^{(p)},B_{r_{p}}}
\end{equation*}
where $\mu ^{(p)}=\mu ^{(p)}-(1,...,r_{p}).$ In particular the computation
of $\varphi _{\ell }(P_{\mu })$ makes appear positive roots corresponding to
a root system of type $B_{r_{p}}.\;$These roots do not belong to the root
lattice associated with $Sp_{2n}$. Hence, there cannot exist an analogue of
Theorem \ref{thC} when $\ell $ is even. With the previous notation, we only
obtain:

\begin{proposition}
Suppose $G=SP_{2n}$ and $\ell =2p.$

\begin{itemize}
\item  If $\mathrm{card}(I^{(0)})\neq \frac{1}{2}\mathrm{card}(J^{(0)}),$ $%
\mathrm{card}(I^{(p)})\neq \frac{1}{2}\mathrm{card}(J^{(p)})$ or there
exists $k\in \{1,...,p-1\}$ such that $\mathrm{card}(X^{(k)})\neq \mathrm{%
card}(J^{(k)})$ then $\varphi _{\ell }(s_{\mu })=0.$

\item  Otherwise, the coefficients appearing in the decomposition of $%
\varphi _{\ell }(s_{\mu })$ on the basis of Weyl characters cannot be
interpreted as branching coefficients and have signs alternatively positive
and negative.
\end{itemize}
\end{proposition}

\subsubsection{For $G=SO_{2n}$}

As for $G=Sp_{2n},$ the coefficients appearing in the decomposition of $%
\varphi_{\ell}(s_{\mu})$ with $\ell=2p$ on the basis of Weyl characters
cannot be interpreted as branching coefficients.\ Note that there is an
additional difficulty in this case. Indeed, $\varphi_{\ell}(P_{\mu})$ cannot
be factorized as a product of polynomials $(1-x^{\beta})$ where $\beta \in%
\mathbb{Z}^{n}$.\ For example, we have for $SO_{4}$ 
\begin{equation*}
\varphi_{2}(P_{(0,0)})=\varphi_{2}\left( (1-\frac{x_{2}}{x_{1}}%
)(1-x_{1}x_{2})\right) =1+x_{2}.
\end{equation*}
This is due to the incompatibility between the signatures defined on the
Weyl groups of types $B$ and $D$ when they are realized as subgroups of the
permutation group $S_{J_{n}}$.

\noindent So we will suppose $\ell=2p-1$ in this paragraph. Recall that the
elements of $W$ are the signed permutations $w$ defined on $J_{n}=\{%
\overline{n},...,\overline{1},1,...,n\}$ such that $\mathrm{card}(\{i\in
I_{n}\mid w(i)<0\})$ is even.\ Set $K_{n}=\{\overline{n-1},...,\overline
{1},0,1,...,n-1\}.\;$Each $w\in W$ can be written $w=\zeta\sigma$ according
to the decomposition of $W$ as the semidirect product $(\mathbb{Z}/2\mathbb{%
Z)}^{n-1}\propto S_{n}$. For any $x\in J_{n}$, we have then $\xi(x)=1$ if $%
w(x)>0$ and $\xi(x)=-1$ otherwise. Given $w\in W$ we define $\widehat
{w}:J_{n}\rightarrow K_{n}$ such that $\widehat{w}(x)=w(x)-\xi(x)$ for any $%
x\in J_{n}.$ Then $\widehat{w}(\overline{x})=\overline{\widehat{w}(x)}$.

\noindent For type $D_{n},$ we have $\rho =\rho _{n}^{\prime
}=(0,1,...,n-1)=\rho _{n}-(1,...,1).\;$Hence 
\begin{equation*}
w\cdot \rho _{n}^{\prime }=w\cdot \rho _{n}-(\xi (1),...,\xi (n)=(\widehat{w}%
(1),...,\widehat{w}(n))=\widehat{w}\cdot \rho _{n}.
\end{equation*}
Then we obtain 
\begin{equation*}
P_{\mu }=x_{1}^{(\mu _{1}+0)}\cdot \cdot \cdot x_{n}^{(\mu
_{n}+n-1)}\sum_{w\in W}\varepsilon (w)x_{1}^{-\widehat{w}(1)}\cdot \cdot
\cdot x_{n}^{-\widehat{w}(n)}.
\end{equation*}
For any $k=0,...,\ell -1$ set 
\begin{equation}
I^{(k)}=(i\in I_{n}\mid \mu _{i}+i-1\equiv k\mathrm{mod}\ell )\text{ and }%
J^{(k)}=(x\in K_{n}\mid x\equiv k\mathrm{mod}\ell ).  \label{IJD}
\end{equation}
We then proceed essentially as in \ref{subsubSP} by using $\widehat{w}$
instead of $w$ and $\rho _{n}^{\prime }=(0,1,...,n-1)$ instead of $\rho
_{n}=(1,...,n).$ We only sketch below the main steps of the computation.

\noindent Set $r_{0}=\mathrm{card}(I^{(0)})$ and for any $k=1,...,p-1$, $%
s_{k}=\mathrm{card}(I_{k}),$ $r_{k}=\mathrm{card}(I_{k})+\mathrm{card}%
(I_{\ell -k})$.\ For $k=1,...,p-1$, $X^{(k)}$ is defined as the increasing
reordering of $\overline{I}^{(k)}\cup I^{(\ell -k)}.$ Consider 
\begin{align*}
\mu ^{(0)}& =\left( \frac{\mu _{i}+i-1}{\ell }\mid i\in I^{(0)}\right) \in 
\mathbb{Z}^{r_{0}}\text{ and for }k>0, \\
\mu ^{(k)}& =\left( \mathrm{sign}(i)\frac{\mu _{\left| i\right| }+\left|
i\right| -1+\mathrm{sign}(i)k}{\ell }\mid i\in X^{(k)}\right) \in \mathbb{Z}%
^{r_{k}}.
\end{align*}
We obtain 
\begin{equation*}
\varphi _{\ell }(P_{\mu })=x_{I^{(0)}}^{\mu
^{(0)}}\prod_{k=1}^{p-1}x_{X^{(k)}}^{\mu ^{(p)}}\times \sum_{w\in
W}\varepsilon (w)\varphi _{\ell }\left( \prod_{i\in I^{(0)}}x_{i}^{-\widehat{%
w}(i)}\prod_{k=1}^{p-1}\prod_{i\in X^{(k)}}x_{i}^{-\widehat{w}(i)-k}\right) .
\end{equation*}
We also have the equivalences 
\begin{equation}
\varphi _{\ell }\left( \prod_{i\in
I^{(0)}}x_{i}^{-w(i)}\prod_{k=1}^{p-1}\prod_{i\in X^{(k)}}^{s_{k}}x_{i}^{-%
\widehat{w}(i)+k}\right) \neq 0\Longleftrightarrow \left\{ 
\begin{tabular}{l}
$\mathrm{(i):}$ $\widehat{w}(I^{(0)}\cup \overline{I}^{(0)})=J^{(0)}$ \\ 
$\mathrm{(ii):}$ $\widehat{w}(X^{(k)})=J^{(\ell -k)}$ for any $k=1,...,p-1$%
\end{tabular}
\right. .  \label{conD}
\end{equation}
We can write $J^{(0)}=(-(r_{0}-1)\ell ,...,0,...,(r_{0}-1)\ell )$ and for $%
k=1,...,p,$ 
\begin{equation*}
J^{(\ell -k)}=(-k-\alpha _{k}\ell ,...,-k+\beta _{k}\ell )\text{, }%
J^{(k)}=(k-\beta _{k}\ell ,...,k+\alpha _{k}\ell )
\end{equation*}
with $\alpha _{k}+\beta _{k}+1=r_{k}.$ Consider $w_{0}\in W$ defined by 
\begin{align}
\widehat{w}_{0}(i_{a}^{(0)})& =(a-1)\ell \text{ for }a\in \{1,...,r_{0}\}
\label{w0D} \\
\widehat{w}_{0}(i_{a}^{(k)})& =-k-\alpha _{k}\ell +(a-1)\ell \text{ for any }%
k=1,...,p-1.  \notag
\end{align}
Denote by $\mathcal{W}$ the set of signed permutations $w\in W$ which verify 
$\mathrm{(i)}$ and $\mathrm{(ii)}$ in (\ref{conD}). We have $w_{0}\in 
\mathcal{W}$.\ Each $w\in \mathcal{W}$ can be written $w=w_{0}v$ where $%
v=(v^{(0)},\tau ^{(1)},...,\tau ^{(p-1)})$ belongs to the direct product $%
W_{I^{(0)}}\times S_{X^{(1)}}\times \cdot \cdot \cdot \times S_{X^{(p-1)}}$
with $W_{I^{(0)}}$ the Weyl group of type $D_{r_{0}}$ defined on $\overline{I%
}^{(0)}\cup I^{(0)}$.\ We have by Lemma \ref{lem_sign} $\varepsilon
(w)=\varepsilon (w_{0})(-1)^{l(v^{(0)})}(-1)^{l(\tau ^{(1)})}\times \cdot
\cdot \cdot \times (-1)^{l(\tau ^{(p-1)})}$.

\noindent Set 
\begin{align*}
P_{0}& =\sum_{v^{(0)}\in W_{I^{(0)}}}(-1)^{l(v^{(0)})}\varphi _{\ell }\left(
\prod_{i\in X^{(0)}}x_{i}^{-\widehat{w}_{0}v^{(0)}(i)}\right) \\
P_{k}& =\sum_{\tau ^{(k)}\in S_{X^{(k)}}}(-1)^{l(\tau ^{(k)})}\varphi _{\ell
}\left( \prod_{i\in X^{(k)}}x_{i}^{-\widehat{w}_{0}v^{(k)}(i)-k}\right) 
\text{ for any }k\in \{1,...,p-1\}.
\end{align*}
We obtain 
\begin{equation*}
\varphi _{\ell }(P_{\mu })=\varepsilon (w_{0})x_{I^{(0)}}^{\mu
^{(0)}}P_{0}\prod_{k=1}^{p-1}x_{X^{(k)}}^{\mu ^{(k)}}P_{k}.
\end{equation*}
Given $v^{(0)}\in W_{I^{(0)}},$ we define $\widehat{v}^{(0)}=v^{(0)}-\xi
_{v} $ where $\xi _{v}(i_{a})=1$ if $v(i_{a})>0$ and $-1$ otherwise.\ By (%
\ref{w0D}), we have for any $a=1,...,r_{0},$ $\widehat{w}%
_{0}v^{(0)}(i_{a}^{(0)})=\widehat{v}^{(0)}(a)\ell $.\ Moreover since $\tau
^{(k)}\in S_{X^{(k)}}$%
\begin{equation*}
\widehat{w}_{0}\tau ^{(k)}(i_{a}^{(k)})=-k-\alpha _{k}\ell +(\tau
^{(k)}(a)-1)\ell \text{ for any }a=1,...,r_{k}.
\end{equation*}
This yields 
\begin{align*}
P_{0}& =\sum_{v^{(0)}\in
W_{I^{(0)}}}(-1)^{l(v^{(0)})}\prod_{a=1}^{r_{0}}x_{i_{a}^{(0)}}^{-\widehat{v}%
^{(0)}(a)}=x_{I^{(0)}}^{-\rho _{r_{0}}^{\prime }}\Delta _{I^{(0)}}\text{ and 
} \\
P_{k}& =\sum_{\tau ^{(k)}\in S_{X^{(k)}}}(-1)^{l(\tau
^{(k)})}\prod_{a=1}^{r_{k}}x_{i_{a}^{(k)}}^{-\tau ^{(k)}(a)+(\alpha
_{k}+1)}=x_{X^{(k)}}^{\eta _{r_{k}}}\Delta _{X^{(k)}}
\end{align*}
where for any $k=1,...,p-1,$ $\eta _{r_{k}}=-\rho _{r_{k}}^{\prime }+(\alpha
_{k},...,\alpha _{k})=-\rho _{r_{k}}+(\alpha _{k}+1,...,\alpha _{k}+1)\in 
\mathbb{Z}^{r_{k}},$%
\begin{align*}
\Delta _{I^{(0)}}& =\prod_{i<j,i,j\in I^{(0)}}(1-\frac{x_{j}}{x_{i}}%
)\prod_{r<s,r,s\in I^{(0)}}(1-x_{r}x_{s})\text{ and} \\
\Delta _{X^{(k)}}& =\prod_{i<j,i,j\in X^{(k)}}(1-\frac{x_{j}}{x_{i}})\text{
for any }k=1,...,p-1.
\end{align*}
This gives 
\begin{equation*}
\varphi _{\ell }(P_{\mu })=\varepsilon (w_{0})x_{I^{(0)}}^{\mu ^{(0)}-\rho
_{r_{0}}^{\prime }}\Delta _{I^{(0)}}\prod_{k=1}^{p-1}x_{X^{(k)}}^{\mu
^{(k)}-\eta _{r_{k}}}\Delta _{X^{(k)}}=\varepsilon (w_{0})x_{I^{(0)}}^{\mu
^{(0)}}\Delta _{I^{(0)}}\prod_{k=1}^{p-1}x_{X^{(k)}}^{\mu ^{(k)}}\Delta
_{X^{(k)}}
\end{equation*}
where 
\begin{equation}
\mu ^{(0)}=\left( \frac{\mu _{i}+i-1}{\ell }\mid i\in I^{(0)}\right)
-(0,...,r_{0}-1)\in \mathbb{Z}^{r_{0}},  \label{mu0D}
\end{equation}
and for any $k=1,...,p-1,$%
\begin{equation}
\mu ^{(k)}=\mid i\in X^{(k)}\left( \mathrm{sign}(i)\frac{\mu _{\left|
i\right| }+\left| i\right| -1+\mathrm{sign}(i)k}{\ell }\mid i\in
X^{(k)}\right) -(0,...,r_{k}-1)+(\alpha _{k},....,\alpha _{k})\in \mathbb{Z}%
^{r_{k}}.  \label{mukD}
\end{equation}
Note that these formulas are essentially the same as for $G=Sp_{2n}$, except
that we use $\rho _{n}^{\prime }=(0,...,n-1)$ instead of $\rho
_{n}=(1,...,n) $ for the half sum of positive roots.\ This gives the
following theorem whose proof is identical to that of Theorem \ref{thC}:

\begin{theorem}
\label{thD}Consider a partition $\mu $ of length $n$ and $\ell =2p-1$ a
positive integer.$\;$Let $I^{(0)}$ and $J^{(0)}$ be as in (\ref{IJD}).\ For
any $k=0,...,p-1$ define the sets $X^{(k)}$ and $J^{(k)}$ by (\ref{deIJC})
and (\ref{IJD}).

\begin{itemize}
\item  If $\mathrm{card}(I^{(0)})\neq \frac{1}{2}(\mathrm{card}(J^{(0)})+1)$
or if there exists $k\in \{1,...,p-1\}$ such that $\mathrm{card}%
(X^{(k)})\neq \mathrm{card}(J^{(k)})$ then $\varphi _{\ell }(s_{\mu })=0.$

\item  Otherwise, set $r_{0}=\mathrm{card}(I^{(0)})$ and for any $%
k=0,...,p-1,$ $r_{k}=\mathrm{card}(X^{(k)})$.\ Let $w_{0}\in W$ be as in (%
\ref{w0D}).\ Consider $\binom{\mu }{\ell }=(\mu ^{(0)},...,\mu ^{(\ell -1)})$
where the $\mu ^{(k)}$'s are defined by (\ref{mu0D}) and (\ref{mukD}). Then $%
\binom{\mu }{\ell }$ is a dominant weight of $P_{\mathcal{I}}^{+}$ with $%
\mathcal{I=}\{I^{(0)},X^{(1)}...,X^{(p-1)}\}$ and we have 
\begin{equation*}
\varphi _{\ell }(s_{\mu })=\varepsilon (w_{0})S_{\binom{\mu }{\ell },%
\mathcal{I}}.
\end{equation*}
\end{itemize}
\end{theorem}

\begin{example}
\ \ \ 

\noindent Consider $\mu =(1,2,3,4,4,4,6,6)$ and take $\ell =3.$ We have $\mu
+\rho _{8}^{\prime }=(1,3,5,7,8,9,12,13).$ Thus $I^{(0)}=\{2,6,7\},X^{(1)}=\{%
\overline{8},\overline{4},\overline{1},3,5\}$ and $J^{(0)}=\{\overline{6},%
\overline{3},0,3,6\},J^{(1)}=\{\overline{5},\overline{2},1,4,7\}$ and $%
J^{(2)}=\{\overline{7},\overline{4},\overline{1},2,5\}.$ In particular $%
\alpha _{1}=2.\;$Then $\mu ^{(0)}=(1,2,2)$ and 
\begin{multline*}
\mu ^{(1)}=\left( -\frac{13-1}{3}-1+3,-\frac{7-1}{3}-2+3,-\frac{1-1}{3}-3+3,%
\frac{5+1}{3}-4+3,\frac{8+1}{3}-5+3\right) \\
=(-2,-1,0,1,1)
\end{multline*}
We have $G_{\mathcal{I}}\simeq SO_{6}\times GL_{5}.$
\end{example}

\subsubsection{For $G=SO_{2n+1}$}

Set $L_{n}=\{\overline{n-1},...,\overline{1},0,1,...,n\}.\;$Each $w\in W$
can be written $w=\zeta\sigma$ according to the decomposition of $W$ as the
semidirect product $(\mathbb{Z}/2\mathbb{Z)}^{n}\propto S_{n}$. Given $w\in
W $ we define $\widetilde{w}:J_{n}\rightarrow L_{n}$ such that $\widetilde
{w}(x)=w(x)+\frac{1}{2}(1-\xi(x))$ for any $x\in J_{n}.$ For any $y\in
L_{n}, $ set $y^{\ast}=\overline{y}+1$.\ We have then $\widetilde{w}%
(\overline {x})=(w(x))^{\ast}=\overline{w(x)}+1$.

\noindent Observe that $\rho =\rho _{n}^{\prime \prime }=(\frac{1}{2},\frac{3%
}{2},...,n-\frac{1}{2})=\rho _{n}-(\frac{1}{2},...,\frac{1}{2})$. Thus 
\begin{equation*}
w\cdot \rho _{n}^{\prime \prime }=w\cdot \rho _{n}-\frac{1}{2}(\xi
(1),...,\xi (n)=(\widetilde{w}(1),...,\widetilde{w}(n))-\frac{1}{2}(1,...,1).
\end{equation*}
This permits to write 
\begin{equation}
P_{\mu }=x_{1}^{(\mu _{1}+1)}\cdot \cdot \cdot x_{n}^{(\mu
_{n}+n)}\sum_{w\in W}\varepsilon (w)x_{1}^{-\widetilde{w}(1)}\cdot \cdot
\cdot x_{n}^{-\widetilde{w}(n)}.  \label{miB}
\end{equation}
For any $k=1,...,\ell $ (observe that $k$ does not run over $\{0,...,\ell
-1\}$ as for $G=Sp_{2n}$ or $SO_{2n}$) set 
\begin{equation}
I^{(k)}=(i\in I_{n}\mid \mu _{i}+i\equiv k\mathrm{mod}\ell )\text{ and }%
J^{(k)}=(x\in L_{n}\mid x\equiv k\mathrm{mod}\ell ).  \label{IJB}
\end{equation}
Note that $(J^{(k)})^{\ast }=J^{(l-k+1)}$. We then proceed essentially as in 
\ref{subsubSP} by using $\widetilde{w}$ instead of $w$.\ We are going to see
that for $G=SO_{2n+1}$, there exists an analogue of Theorem \ref{thC}
whatever the parity of $\ell $.

\paragraph{The even case $\ell=2p$}

For any $k=1,...,p$, set $s_{k}=\mathrm{card}(I^{(k)}),$ $r_{k}=\mathrm{card}%
(I^{(k)})+\mathrm{card}(I^{(\ell -k+1)})$ and define $X^{(k)}$ as the
increasing reordering of $\overline{I}^{(k)}\cup I^{(\ell -k+1)}.\;$Set 
\begin{equation}
X^{(k)}=(i_{1}^{(k)},...,i_{r_{k}}^{(k)}).  \label{defXKB}
\end{equation}
For $k=1,...,p$ consider the $r_{k}$-tuple $\mu ^{(k)}\ $such that 
\begin{equation*}
\mu ^{(k)}=\left( \mathrm{sign}(i)\frac{\mu _{\left| i\right| }+\left|
i\right| +\mathrm{sign}(i)k-\frac{1+\mathrm{sign}(i)}{2}}{\ell }\mid i\in
X^{(k)}\right) \in \mathbb{Z}^{r_{k}}.
\end{equation*}
For any $i\in I^{(k)}$ with $k=1,...,p,$ we have $x_{i}^{-\widetilde{w}%
(i)-1}=x_{\overline{i}}^{-\widetilde{w}(\overline{i})}.$ Thus 
\begin{equation*}
\prod_{i\in X^{(k)}}x_{i}^{-\widetilde{w}(i)}=\prod_{i\in \overline{I}%
^{(k)}}x_{i}^{-\widetilde{w}(i)}\prod_{i\in I^{(\ell -k+1)}}x_{i}^{-%
\widetilde{w}(i)}=\prod_{i\in I^{(k)}}x_{i}^{-\widetilde{w}(i)-1}\prod_{i\in
I^{(\ell -k+1)}}x_{i}^{-\widetilde{w}(i)}
\end{equation*}
and by definition of the $\mu ^{(k)}$'s, (\ref{miB}) can be rewritten\ 
\begin{equation*}
P_{\mu }=\prod_{k=1}^{p}x_{X^{(k)}}^{\ell \mu ^{(k)}}\times \sum_{w\in
W}\varepsilon (w)\prod_{k=1}^{p}\prod_{i\in X^{(k)}}x_{i}^{-\widetilde{w}%
(i)}\times \prod_{k=1}^{p}\prod_{i\in X^{(k)}}x_{i}^{-k+1}.
\end{equation*}
We obtain 
\begin{equation*}
\varphi _{\ell }(P_{\mu })=\prod_{k=1}^{p}x_{X^{(k)}}^{\mu ^{(k)}}\times
\sum_{w\in W}\varepsilon (w)\varphi _{\ell }\left(
\prod_{k=1}^{p}\prod_{i\in X^{(k)}}x_{i}^{-\widetilde{w}(i)-k+1}\right)
\end{equation*}
We deduce the equivalences 
\begin{equation}
\varphi _{\ell }\left( \prod_{k=0}^{p}\prod_{i\in
X^{(k)}}^{s_{k}}x_{i}^{-w(i)+k-1}\right) \neq 0\Longleftrightarrow 
\widetilde{w}(X^{(k)})=J^{(\ell -k+1)}\text{ for any }k=1,...,p.
\label{conB}
\end{equation}
In particular we must have $\mathrm{card}(J^{(\ell -k+1)})=\mathrm{card}%
(J^{(k)})=r_{k}$. We can write 
\begin{equation*}
J^{(\ell -k+1)}=(-k+1-\alpha _{k}\ell ,...,-k+1+\beta _{k}\ell )\text{ and }%
J^{(k)}=(k-\beta _{k}\ell ,...,k+\alpha _{k}\ell )
\end{equation*}
with $\alpha _{k}+\beta _{k}+1=r_{k}.$ Consider $w_{0}\in W$ defined by 
\begin{equation}
\widetilde{w}_{0}(i_{a}^{(k)})=-k+1-\alpha _{k}\ell +(a-1)\ell \text{ for
any }k=1,...,p.  \label{w0B}
\end{equation}
Denote by $\mathcal{W}$ the set of signed permutations $w\in W$ which verify
the right-hand side of (\ref{conB}). We have $w_{0}\in \mathcal{W}$.\ Each $%
w\in \mathcal{W}$ can be written $w=w_{0}v$ where $\tau =(\tau
^{(1)},...,\tau ^{(p)})$ belongs to the direct product $S_{X^{(1)}}\times
\cdot \cdot \cdot \times S_{X^{(p)}}$.\ We have also by Lemma \ref{lem_sign} 
$\varepsilon (w)=\varepsilon (w_{0})(-1)^{l(\tau ^{(1)})}\times \cdot \cdot
\cdot \times (-1)^{l(\tau ^{(p)})}$.

\noindent For any $k=1,...,p,$ set 
\begin{equation*}
P_{k}=\sum_{\tau ^{(k)}\in S_{X^{(k)}}}(-1)^{l(\tau ^{(k)})}\varphi _{\ell
}\left( \prod_{i\in X^{(k)}}x_{i}^{-\widetilde{w}_{0}\tau
^{(k)}(i)-k+1}\right) .
\end{equation*}
We obtain 
\begin{equation*}
\varphi _{\ell }(P_{\mu })=\varepsilon
(w_{0})\prod_{k=1}^{p}x_{X^{(k)}}^{\mu ^{(k)}}P_{k}.
\end{equation*}
By (\ref{w0B}) we have 
\begin{equation*}
w_{0}\tau ^{(k)}(i_{a}^{(k)})=-k+1-\alpha _{k}\ell +(\tau ^{(k)}(a)-1)\ell 
\text{ for any }a=1,...,r_{k}.
\end{equation*}
This yields 
\begin{equation*}
P_{k}=\sum_{\tau ^{(k)}\in S_{X^{(k)}}}(-1)^{l(\tau
^{(k)})}\prod_{a=1}^{r_{k}}x_{i_{a}^{(k)}}^{-\tau ^{(k)}(a)+(\alpha
_{k}+1)}=x_{X^{(k)}}^{\eta _{r_{k}}}\Delta _{X^{(k)}}
\end{equation*}
where for any $k=1,...,p,$ $\eta _{r_{k}}=-\rho _{r_{k}}+(\alpha
_{k}+1,...,\alpha _{k}+1)\in \mathbb{Z}^{r_{k}}$ and 
\begin{equation*}
\Delta _{X^{(k)}}=\prod_{i<j\ i,j\in X^{(k)}}(1-\frac{x_{j}}{x_{i}}).
\end{equation*}
Note that the computation only makes root systems of type $A$ appear in this
case. This gives 
\begin{equation*}
\varphi _{\ell }(P_{\mu })=\varepsilon
(w_{0})\prod_{k=1}^{p}x_{X^{(k)}}^{\mu ^{(k)}-\eta _{r_{k}}}\Delta
_{X^{(k)}}=\varepsilon (w_{0})\prod_{k=1}^{p}x_{X^{(k)}}^{\mu ^{(k)}}\Delta
_{X^{(k)}}
\end{equation*}
where for any $k=1,...,p,$%
\begin{equation}
\mu ^{(k)}=\left( \mathrm{sign}(i)\frac{\mu _{\left| i\right| }+\left|
i\right| +\mathrm{sign}(i)k-\frac{1+\mathrm{sign}(i)}{2}}{\ell }\mid i\in
X^{(k)}\right) -(1,...,r_{k})+(\alpha _{k+1},....,\alpha _{k+1})\in \mathbb{Z%
}^{r_{k}}.  \label{mukB}
\end{equation}
Similarly to Theorem \ref{thC} we obtain:

\begin{theorem}
\label{thB1}Consider a partition $\mu $ of length $n$ and $\ell =2p$ a
positive integer.$\;$For any $k=1,...,p$ define the sets $X^{(k)},J^{(k)}$
by (\ref{IJB}) and (\ref{defXKB}).

\begin{itemize}
\item  If there exists $k\in \{1,...,p\}$ such that $\mathrm{card}%
(X^{(k)})\neq \mathrm{card}(J^{(k)})$ then $\varphi _{\ell }(s_{\mu })=0.$

\item  Otherwise, for any $k=1,...,p,$ set $r_{k}=\mathrm{card}(X^{(k)})$.\
Let $w_{0}\in W$ be as in (\ref{w0B}).\ Consider $\binom{\mu }{\ell }=(\mu
^{(1)},...,\mu ^{(p)})$ where the $\mu ^{(k)}$'s are defined by (\ref{mukB}%
). Then $\binom{\mu }{\ell }$ is a dominant weight of $P_{\mathcal{I}}^{+}$
with $\mathcal{I=}\{X^{(1)}...,X^{(p)}\}$ and we have 
\begin{equation*}
\varphi _{\ell }(s_{\mu })=\varepsilon (w_{0})S_{\binom{\mu }{\ell },%
\mathcal{I}}.
\end{equation*}
\end{itemize}
\end{theorem}

\begin{example}
Consider $\mu =(2,5,5,6,7,9)$ and $\ell =2.\;$Then $\mu +\rho
_{6}=(3,7,8,10,12,15).\;$Hence $I_{1}=\{1,2,6\}$ and $I_{2}=\{3,4,5\}.\;$%
Moreover $J_{2}=\{\overline{4},\overline{2},0,2,4,6\}$ and $J_{1}=\{%
\overline{5},\overline{3},\overline{1},1,3,5\}.$ Then $\widetilde{w}_{0}$
sends $X_{1}=\{\overline{6},\overline{2},\overline{1},3,4,5\}$ on $J_{2}.$
This gives 
\begin{equation*}
\widetilde{w}_{0}=\left( 
\begin{array}{cccccccccccc}
\overline{6} & \overline{5} & \overline{4} & \overline{3} & \overline{2} & 
\overline{1} & 1 & 2 & 3 & 4 & 5 & 6 \\ 
\overline{4} & \overline{5} & \overline{3} & \overline{1} & \overline{2} & 0
& 1 & 3 & 2 & 4 & 6 & 5
\end{array}
\right)
\end{equation*}
by using (\ref{rel}).\ Hence 
\begin{equation*}
w_{0}=\left( 
\begin{array}{cccccccccccc}
\overline{6} & \overline{5} & \overline{4} & \overline{3} & \overline{2} & 
\overline{1} & 1 & 2 & 3 & 4 & 5 & 6 \\ 
\overline{5} & \overline{6} & \overline{4} & \overline{2} & \overline{3} & 
\overline{1} & 1 & 3 & 2 & 4 & 6 & 5
\end{array}
\right) .
\end{equation*}
We have $\varepsilon (\mu )=1,$ $\alpha _{1}=2$ and 
\begin{equation*}
\mu ^{(1)}=(-7,-3,-1,4,5,6)-(1,2,3,4,5,6)+(3,3,3,3,3,3)=(-5,-2,-1,3,3,3).
\end{equation*}
We have then $G_{\mathcal{I}}\simeq GL_{6}.$
\end{example}

\paragraph{The case $\ell=2p+1$}

In addition to the sets $X^{(k)},k=1,...,p$ defined in (\ref{defXKB}), we
have also to consider $I^{(p+1)}=\{i_{1}^{(p+1)},...,i_{r_{p+1}}^{(p+1)}\}.$
This yields to define 
\begin{equation*}
\mu ^{(p+1)}=\left( \frac{\mu _{i}+i+p}{\ell }\mid i\in I^{(p+1)}\right)
\end{equation*}
We have\ 
\begin{equation*}
\varphi _{\ell }(P_{\mu })=x_{I^{(p+1)}}^{\mu
^{(p+1)}}\prod_{k=1}^{p}x_{X^{(k)}}^{\mu ^{(k)}}\times \sum_{w\in
W}\varepsilon (w)\varphi _{\ell }\left( \prod_{i\in I^{(p+1)}}x_{i}^{-%
\widetilde{w}(i)-p}\prod_{k=1}^{p}\prod_{i\in X^{(k)}}x_{i}^{-\widetilde{w}%
(i)-k-1}\right)
\end{equation*}
and the equivalence 
\begin{equation}
\varphi _{\ell }\left( \prod_{i\in I^{(p+1)}}x_{i}^{-\widetilde{w}%
(i)-p}\prod_{k=1}^{p}\prod_{i\in X^{(k)}}x_{i}^{-w(i)+k-1}\right) \neq
0\Longleftrightarrow \left\{ 
\begin{tabular}{l}
$\widetilde{w}(X^{(k)})=J^{(\ell -k+1)}\text{ for any }k=1,...,p$ \\ 
$\widetilde{w}(I^{(p+1)}\cup \overline{I}^{(p+1)})=J^{(\ell -p)}=J^{(p+1)}$%
\end{tabular}
\right. .  \label{condB}
\end{equation}
Indeed $(J^{(p+1)})^{\ast }=J^{(p+1)}$.\ In particular we must have $\mathrm{%
card}(J^{(p+1)})=2\mathrm{card}(I^{(p+1)})=2r_{p+1}$. Thus we can set $%
J^{(p+1)}=(-p-(r_{p+1}-1)\ell ,...,-p+r_{p+1}\ell ).$ Consider $w_{0}\in W$
defined by (\ref{w0B}) and 
\begin{equation}
\widetilde{w}_{0}(i_{a}^{(p+1)})=-p+a\ell \text{ for any }a=1,...,r_{p+1}.
\label{w0B2}
\end{equation}
Denote by $\mathcal{W}$ the set of signed permutations $w\in W$ which verify
the right-hand side of (\ref{condB}). We have $w_{0}\in \mathcal{W}$.\ Each $%
w\in \mathcal{W}$ can be written $w=w_{0}v$ where $v=(\tau ^{(1)},...,\tau
^{(p)},v^{(p+1)})$ belongs to the direct product $S_{X^{(1)}}\times \cdot
\cdot \cdot \times S_{X^{(p)}}\times W_{I^{(p+1)}}$.\ We have also $%
\varepsilon (w)=\varepsilon (w_{0})(-1)^{l(\tau ^{(1)})}\times \cdot \cdot
\cdot \times (-1)^{l(\tau ^{(p)})}(-1)^{l(v^{(p+1)})}$. This permits to
write 
\begin{align*}
\varphi _{\ell }(P_{\mu })& =\varepsilon (w_{0})x_{I^{(p+1)}}^{\mu
^{(p+1)}}P_{p+1}\prod_{k=1}^{p}x_{X^{(k)}}^{\mu ^{(k)}}P_{k}\text{ where} \\
P_{p+1}& =\sum_{v^{(p+1)}\in W_{I^{(p+1)}}}(-1)^{l(v^{(p+1)})}\varphi _{\ell
}\left( \prod_{i\in I^{(p+1)}}x_{i}^{-\widetilde{w}_{0}v^{(p+1)}(i)-p}%
\right) .
\end{align*}
The functions $P_{k},k=1,...,p$ can be computed as in the even case.\ For $%
P_{p+1}$, observe that each $v^{(p+1)}\in W_{I^{(p+1)}}$ can be written $%
v^{(p+1)}=\zeta \sigma $ with $\sigma \in S_{I^{(p+1)}}$. According to this
decomposition we have for any $a=1,...,r_{p+1}$, $\widetilde{w}%
_{0}v^{(p+1)}(i_{a}^{(p+1)})=\xi (a)(-p+\sigma (a)\ell ).$%
\begin{multline*}
P_{p+1}=\sum_{v^{(p+1)}\in W_{I^{(p+1)}}}(-1)^{l(v^{(p+1)})}\varphi _{\ell
}\left( \prod_{a=1}^{r_{p+1}}x_{i_{a}^{(p+1)}}^{-\xi (a)(-p+\sigma (a)\ell
)-p}\right) = \\
\sum_{v^{(p+1)}\in
W_{I^{(p+1)}}}(-1)^{l(v^{(p+1)})}\prod_{a=1}^{r_{p+1}}x_{i_{a}^{(p+1)}}^{-%
\frac{1-\xi (a)}{2}-\xi \sigma (a)}.
\end{multline*}
Thus 
\begin{equation*}
P_{p+1}=\prod_{i\in I^{(p+1)}}x_{i}^{-1/2}\sum_{v^{(p+1)}\in
W_{I^{(p+1)}}}(-1)^{l(v^{(p+1)})}\left( \nu ^{(p+1)}\cdot
\prod_{a=1}^{r_{p+1}}x_{i_{a}^{(p+1)}}^{-(a-\frac{1}{2})}\right)
=x_{I^{(p+1)}}^{-\rho _{r_{p+1}}}\Delta _{I^{(p+1)}}
\end{equation*}
\begin{equation*}
\Delta _{I^{(p+1)}}=\prod_{i<j\ i,j\in I^{(p+1)}}(1-\frac{x_{j}}{x_{i}}%
)\prod_{r<s\ r,s\in I^{(p+1)}}(1-x_{r}x_{s})\prod_{i\in I^{(p+1)}}(1-x_{i}).
\end{equation*}
This means that when $\ell $ is odd 
\begin{equation*}
\varphi _{\ell }(P_{\mu })=\varepsilon
(w_{0})\prod_{k=1}^{p}x_{X^{(k)}}^{\mu ^{(k)}}\Delta _{X^{(k)}}\times
x_{I^{(p+1)}}^{\mu ^{(p+1)}}\Delta _{I^{(p+1)}}
\end{equation*}
where 
\begin{equation}
\mu ^{(p+1)}=\left( \frac{\mu _{i}+i+p}{\ell }\mid i\in I^{(p+1)}\right)
-(1,...,r_{p+1})\in \mathbb{Z}^{r_{p+1}}.  \label{mipB}
\end{equation}
This gives the following theorem:

\begin{theorem}
\label{thB2}Consider a partition $\mu $ of length $n$ and $\ell =2p+1$ a
positive integer.$\;$Define $X^{(k)},J^{(k)}$ $k=1,...,p$ and $%
I^{(p+1)},J^{(p+1)}$ by (\ref{IJB}) and (\ref{defXKB}).

\begin{itemize}
\item  If $\mathrm{card}(I^{(p+1)})\neq \frac{1}{2}\mathrm{card}(J^{(p+1)})$
or if there exists $k\in \{1,...,p\}$ such that $\mathrm{card}(X^{(k)})\neq 
\mathrm{card}(J^{(k)})$ then $\varphi _{\ell }(s_{\mu })=0.$

\item  Otherwise, set $r_{p+1}=\mathrm{card}(I^{(p+1)})$ and for any $%
k=1,...,p,$ $r_{k}=\mathrm{card}(X^{(k)})$.\ Let $w_{0}\in W$ verifying (\ref
{w0B}) and (\ref{w0B2}).\ Consider $\binom{\mu }{\ell }=(\mu ^{(p+1)},\mu
^{(1)},...,\mu ^{(p)})$ where the $\mu ^{(k)}$'s are defined by (\ref{mukB})
and (\ref{mipB}). Then $\binom{\mu }{\ell }$ is a dominant weight of $P_{%
\mathcal{I}}^{+}$ with $\mathcal{I=}\{I^{(p+1)},X^{(1)},...,X^{(p)}\}$ and
we have 
\begin{equation*}
\varphi _{\ell }(s_{\mu })=\varepsilon (w_{0})S_{\binom{\mu }{\ell },%
\mathcal{I}}.
\end{equation*}
\end{itemize}
\end{theorem}

\begin{example}
Consider $\mu =(1,5,5,6,7,9)$ and take $\ell =3.$ We have $\mu +\rho
_{6}=(2,7,8,10,12,15).$ Thus $X^{(1)}=\{\overline{4},\overline{2}%
,5,6\},I^{(2)}=\{1,3\}$ and $J^{(1)}=\{\overline{5},\overline{2}%
,1,4\},J^{(2)}=\{\overline{4},\overline{1},2,5\}$. In particular $\alpha
_{2}=1.\;$Then 
\begin{equation*}
\mu ^{(1)}=\left( -\frac{10-1}{3}-1+2,-\frac{7-1}{3}-2+2,\frac{12}{3}-3+2,%
\frac{15}{3}-4+2\right) =(-2,-2,3,3)
\end{equation*}
and $\mu ^{(2)}=(\frac{2+1}{3}-1,\frac{8+1}{3}-2)=(0,1).$ Moreover, one has
by using (\ref{rel}) 
\begin{equation*}
\widetilde{w}_{0}=\left( 
\begin{array}{cccccccccccc}
\overline{6} & \overline{5} & \overline{4} & \overline{3} & \overline{2} & 
\overline{1} & 1 & 2 & 3 & 4 & 5 & 6 \\ 
\overline{3} & 0 & \overline{5} & \overline{4} & \overline{2} & \overline{1}
& 2 & 3 & 5 & 6 & 1 & 4
\end{array}
\right) .
\end{equation*}
Hence 
\begin{equation*}
w_{0}=\left( 
\begin{array}{cccccccccccc}
\overline{6} & \overline{5} & \overline{4} & \overline{3} & \overline{2} & 
\overline{1} & 1 & 2 & 3 & 4 & 5 & 6 \\ 
\overline{4} & \overline{1} & \overline{6} & \overline{5} & \overline{3} & 
\overline{2} & 2 & 3 & 5 & 6 & 1 & 4
\end{array}
\right)
\end{equation*}
and $\varepsilon (\mu )=1.$ We have moreover $G_{\mathcal{I}}\simeq
SO_{5}\times GL_{4}.$
\end{example}

\section{Parabolic Kazhdan-Lusztig polynomials}

We recall briefly in this section some basics on Affine Hecke algebras and
parabolic Kazhdan-Lusztig polynomials associated with classical root
systems.\ The reader is referred to \cite{NR} and \cite{So} for detailed
expositions. Note that the definition of the Hecke algebra used in \cite{NR}
coincides with that used in \cite{LT} and \cite{So} (with generators $H_{w}$%
) up to the change $q\rightarrow q^{-1}$.

\subsection{Extended affine Weyl group\label{exWG}}

Consider a root system of type $A_{n-1},B_{n},C_{n}$ or $D_{n}$. For any $%
\beta \in P\cap \mathbb{Z}^{n}$, we denote by $t_{\beta }$ the translation
defined in $\frak{h}_{\mathbb{R}}^{\ast }$ by $\gamma \longmapsto \gamma
+\beta $. The extended affine Weyl group $\widehat{W}$ is the group 
\begin{equation*}
\widehat{W}=\{wt_{\beta }\mid w\in W,\beta \in P\}
\end{equation*}
with multiplication determined by the relations $t_{\beta }t_{\gamma
}=t_{\beta +\gamma }$ and $wt_{\beta }=t_{w\cdot \beta }w.$ The group $%
\widehat{W}$ is not a Coxeter group but contains the affine Weyl group $%
\widetilde{W}$ generated by reflections through the affine hyperplanes $%
H_{\alpha ,k}=\{\beta \in \frak{h}_{\mathbb{R}}^{\ast }\mid (\beta ,\alpha
^{\vee })=k\}.\;$It makes sense to define a length function on $\widehat{W}$
verifying 
\begin{equation}
l(wt_{\beta })=\sum_{\alpha \in R_{+}}\left| (\beta ,\alpha ^{\vee
})+1_{R_{-}}(w\cdot \alpha )\right|  \label{defl}
\end{equation}
where for any $w\in W$, $1_{R_{-}}(w\cdot \alpha )=0$ if $w\cdot \alpha \in
R_{+}$ and $1_{R_{-}}(w\cdot \alpha )=1$ if $w\cdot \alpha \in -R_{+}=R_{-}$%
. Write $n_{\beta }$ for the element of maximal length in $Wt_{\beta }W$. It
follows from (\ref{defl}) that for any $\lambda \in P_{+}$ we have $%
l(wt_{\lambda })=l(w)+l(t_{\lambda }).$ This gives 
\begin{equation}
n_{\lambda }=w_{0}t_{\lambda }  \label{nlmabda}
\end{equation}
where $w_{0}$ denotes the longest element of $W$.\ There exists a unique
element $\eta \in R_{+}$ such that the fundamental alcove 
\begin{equation*}
\mathcal{\mathcal{A}}=\{\beta \in \frak{h}_{\mathbb{R}}^{\ast }\mid (\beta
,\alpha ^{\vee })\geq 0\ \forall \alpha \in R_{+}\text{ and }(\beta ,\eta
^{\vee })<1\}
\end{equation*}
is a fundamental region for the action of $\widetilde{W}$ on $\frak{h}_{%
\mathbb{R}}^{\ast }$.\ This means that, for any $\beta \in \frak{h}_{\mathbb{%
R}}^{\ast }$, the orbit $\widetilde{W}\cdot \beta $ intersects $\mathcal{A}$
in a unique point. Each $w\in \widehat{W}$ can be written on the form $w=w_{%
\mathcal{A}}w_{\mathrm{aff}}$ where $w_{\mathrm{aff}}\in \widetilde{W}$ and $%
w_{\mathcal{A}}$ belongs to the stabilizer of $\mathcal{A}$ under the action
of $\widehat{W}$.\ This implies that $\mathcal{A}$ is also a fundamental
domain for the action of $\widehat{W}$ on $\frak{h}_{\mathbb{R}}^{\ast }$.\
The Bruhat ordering on $\widehat{W}$ is defined by taking the transitive
closure of the relations 
\begin{equation*}
w<sw\text{ whenever }l(w)<l(sw)
\end{equation*}
for all $w\in \widehat{W}$ and all (affine) reflections $s\in \widetilde{W}$.

\noindent In fact the natural action of $\widehat{W}$ on the weight lattice $%
P$ obtained by considering $P$ as a sublattice of $\frak{h}_{\mathbb{R}%
}^{\ast}$ is not that which is relevant for our purpose. For any integer $%
m\in\mathbb{Z}^{\ast}$ we obtain a faithful representation $\pi_{m}$ of $%
\widehat{W}$ on $P$ by setting for any $\beta,\gamma\in P,w\in W$%
\begin{equation*}
\pi_{m}(w)\cdot\gamma=w\cdot\gamma\text{ and }\pi_{m}(t_{\beta})\cdot
\gamma=\gamma+m\beta.
\end{equation*}
\textbf{Warning:} In the sequel, the extended affine Weyl group $\widehat{W}$
acts on the weight lattice $P$ via $\pi_{-\ell}$ where $\ell$ is a fixed
nonnegative integer.

\bigskip

\noindent We write for simplicity $wt_{\beta }\cdot \gamma $ rather that $%
\pi _{-\ell }(wt_{\beta })\cdot \gamma $.$\;$Hence for any $w\in W$ and any $%
\beta \in P,$ we have $wt_{\beta }\cdot \gamma =w\cdot \gamma -\ell w\cdot
\beta $.\ The fundamental region for this new action of $\widehat{W}$ on $P$
is the alcove $\mathcal{A}_{\ell }$ obtained by expanding $\mathcal{A}$ with
the factor $-\ell .$ This gives 
\begin{equation*}
\mathcal{A}_{\ell }=\left\{ 
\begin{array}{l}
\{\nu =(\nu _{1},...,\nu _{n})\mid 0\geq \nu _{1}\geq \cdot \cdot \cdot \geq
\nu _{n}>-\ell \}\text{ for types }A \\ 
\{\nu =(\nu _{1},...,\nu _{n})\mid 0\geq \nu _{1}\geq \cdot \cdot \cdot \geq
\nu _{n}\geq -\ell /2\}\text{ for types }B,C \\ 
\{\nu =(\nu _{1},...,\nu _{n})\mid 0\geq \nu _{1}\geq \nu _{2}\geq \cdot
\cdot \cdot \geq \nu _{n}\geq -\ell /2\}\sqcup \\ 
\{\nu =(\nu _{1},...,\nu _{n})\mid 0>\nu _{1}\geq \nu _{2}\geq \cdot \cdot
\cdot \geq \nu _{n}>-\ell /2\}\text{ for type }D
\end{array}
\right. .
\end{equation*}
Consider a weight $\beta \in P$.\ Then its orbit intersects $\mathcal{A}%
_{\ell }$ in a unique weight $\nu .$ Then there is a unique $w(\beta )\in 
\widehat{W}$ of minimal length such that $w(\beta )\cdot \nu =\beta $. We
denote by $W_{\nu }$ the stabilizer of $\nu \in \mathcal{A}_{\ell }$ in $%
\widehat{W}$. Since $\nu \in \mathcal{A}_{\ell }$, $W_{\nu }$ is in fact a
subgroup of $W$.\ 

\begin{lemma}
\label{Lem_nlambda}Consider $\lambda \in P^{+}$ and suppose $\ell >n$. Then

\begin{enumerate}
\item  $w(\ell \lambda +\rho )=n_{\lambda ^{\ast }}\tau ^{-n+1}$ with $%
\lambda ^{\ast }=-w_{0}(\lambda )$ and $\tau =s_{1}s_{2}\cdot \cdot \cdot
s_{r-1}t_{\varepsilon _{1}}$ for type $A.$

\item  $w(\ell \lambda +\rho )=n_{\lambda }$ for types $B,C$ and $D.$
\end{enumerate}
\end{lemma}

\begin{proof}
$1:$ See Lemma 2.3 in \cite{LT}.

$2:$ Observe first that $w_{0}\cdot \rho =-\rho $ belongs to $\mathcal{A}%
_{\ell }$ for types $B,C,D$ since $\ell >n.$ We have 
\begin{equation*}
\ell \lambda +\rho =t_{-\lambda }\cdot \rho =t_{-\lambda }w_{0}\cdot
(w_{0}\cdot \rho )=t_{-\lambda }w_{0}\cdot (-\rho ).
\end{equation*}
Moreover $W_{-\rho }=\{1\}.$ Since $-\rho \in \mathcal{A}_{\ell }$ this
means that $w(\ell \lambda +\rho )=t_{-\lambda }w_{0}=w_{0}t_{w_{0}\cdot
(-\lambda )}=w_{0}t_{\lambda }=n_{\lambda }$ where the last equality follows
from (\ref{nlmabda}).
\end{proof}

\subsection{Affine Hecke algebra and K-L polynomials}

The Hecke algebra associated with the root system $R$ of type $%
A_{n},B_{n},C_{n}$ or $D_{n}$ is the $\mathbb{Z}[q,q^{-1}]$-algebra defined
by the generators $T_{w},w\in \widehat{W}$ and relations 
\begin{align*}
T_{w_{1}}T_{w_{2}}& =T_{w_{1}}T_{w_{2}}\text{ if }%
l(w_{1}w_{2})=l(w_{1})+l(w_{2}), \\
T_{s_{i}}T_{w}& =(q^{-1}-q)T_{w}+T_{s_{i}w}\text{ if }l(s_{i}w)<l(w)\text{
and }i\in I_{n}.
\end{align*}
In particular we have $T_{i}^{2}=(q^{-1}-q)T_{i}+1$ for any $i\in I_{n}.$
The bar involution on $\widehat{H}$ is the $\mathbb{Z}$-linear automorphism
defined by 
\begin{equation*}
\overline{q}=q^{-1}\text{ and }\overline{T}_{w}=T_{w^{-1}}^{-1}\text{ for
any }w\in \widehat{W}.
\end{equation*}
Kazhdan and Lusztig have proved that there exists a unique basis $%
\{C_{w}^{\prime }\mid w\in \widehat{W}\}$ of $\widehat{H}$ such that 
\begin{equation*}
\overline{C}_{w}^{\prime }=C_{w}^{\prime }\text{ and }C_{w}^{\prime
}=\sum_{y\leq w}p_{y,w}T_{y}
\end{equation*}
where $p_{w,w}=1$ and $p_{y,w}\in q\mathbb{Z}[q]$ for any $y<w.$ We will
refer to the polynomials $p_{y,w}(q)$ as Kazhdan-Lusztig polynomials. They
are renormalizations of the polynomials $P_{y,w}$ originally introduced by
Kazhdan and Lusztig in \cite{KL}. Namely we have $%
p_{y,w}=q^{l(w)-l(y)}P_{y,w}.$

\noindent Let us define the $q$-partition function $\mathcal{P}_{q}$ by 
\begin{equation*}
\prod_{\alpha \in R_{+}}\frac{1}{1-qx^{\alpha }}=\sum_{\beta \in \mathbb{Z}%
^{n}}\mathcal{P}_{q}(\beta )x^{\beta }.
\end{equation*}
Given $\lambda $ and $\mu $ in $P$, the Lusztig $q$-analogue $K_{\lambda
,\mu }(q)$ is defined by 
\begin{equation*}
K_{\lambda ,\mu }(q)=\sum_{w\in W}\varepsilon (w)\mathcal{P}_{q}(w\circ
\lambda -\mu ).
\end{equation*}
Then one has the following theorem due to Lusztig:

\begin{theorem}
\label{THL}Suppose $\lambda ,\mu $ are dominant weights. Then $K_{\lambda
,\mu }(q)=p_{n_{\mu },n_{\mu }}(q)$.
\end{theorem}

\bigskip

One defines the action of the bar involution on the parabolic module $P_{\nu
}=\widehat{H}\nu $, $\nu \in \mathcal{A}_{\ell }$, by setting $\overline{q}%
=q^{-1}$ and $\overline{h\cdot \nu }=\overline{h}\cdot \nu $ for any $h\in 
\widehat{H}.\;$Deodhar has proved that there exist two bases $\{C_{\lambda
}^{+}\mid \lambda \in \widehat{W}\cdot \nu \}$ and $\{C_{\lambda }^{-}\mid
\lambda \in \widehat{W}\cdot \nu \}$ of $P_{\nu }$ belonging respectively to 
\begin{equation*}
L_{\nu }^{+}=\bigoplus_{\lambda \in \widehat{W}\cdot \nu }\mathbb{Z}%
[q]\lambda \text{ and }L_{\nu }^{-}=\bigoplus_{\lambda \in \widehat{W}\cdot
\nu }\mathbb{Z}[q^{-1}]\lambda
\end{equation*}
characterized by 
\begin{equation*}
\left\{ 
\begin{array}{l}
\overline{C}_{\lambda }^{+}=C_{\lambda }^{+} \\ 
C_{\lambda }^{+}\equiv \lambda \mathrm{mod}qL_{\nu }^{+}
\end{array}
\right. \text{ and }\left\{ 
\begin{array}{l}
\overline{C}_{\lambda }^{-}=C_{\lambda }^{-} \\ 
C_{\lambda }^{-}\equiv \lambda \mathrm{mod}q^{-1}L_{\nu }^{-}
\end{array}
\right. .
\end{equation*}
We will only need the basis $\{C_{\lambda }^{-}\mid \lambda \in \widehat{W}%
\cdot \nu \}$ in the sequel. The parabolic Kazhdan-Lusztig polynomials $%
P_{\lambda ,\mu }^{-}$ are then defined by the expansion 
\begin{equation*}
C_{\lambda }^{-}=\sum_{\mu \in \widehat{W}\cdot \lambda }(-1)^{l(w(\lambda
))+l(w(\mu ))}P_{\mu ,\lambda }^{-}(q^{-1})\mu
\end{equation*}
(see \cite{So} Theorem 3.5). In particular they belong to $\mathbb{Z}[q]$.\
Their expansion in terms of the ordinary Kazhdan-Lusztig polynomials is
given by the following theorem due to Deodhar:

\begin{theorem}
\label{th_decdeo}Consider $\nu \in \mathcal{A}_{\ell }$ and $\lambda \in 
\widehat{W}\cdot \nu .$ Then for any $\mu \in \widehat{W}\cdot \lambda $ we
have 
\begin{equation}
P_{\lambda ,\mu }^{-}(q)=\sum_{z\in W_{\nu }}(-q)^{l(z)}p_{w(\mu
)z,w(\lambda )}(q)  \label{decdeo}
\end{equation}
with the notation of \ref{exWG}.
\end{theorem}

\noindent\textbf{Remark: }When $\nu$ is regular, that is $W_{\nu}=\{1\},$ we
have $P_{\lambda,\mu}^{-}(q)=p_{w(\mu),w(\lambda)}(q).$

\section{Generalized Hall-Littlewood functions}

\subsection{Plethysm and parabolic K-L polynomials}

Consider $\ell$ a nonnegative integer and $\zeta\in\mathbb{C}$ such that $%
\zeta^{2}$ is a primitive $\ell$-th root of $1$. We briefly recall in this
paragraph the arguments of \cite{LT} which establish that the coefficients
of the plethysm $\psi_{\ell}(s_{\lambda})$ on the basis of Weyl characters
are, up to a sign, parabolic Kazhdan-Lusztig polynomials specialized at $q=1$%
.

\noindent For any $\lambda\in P_{+},$ denote by $V_{q}(\lambda)$ the finite
dimensional $U_{q}(\frak{g})$-module of highest weight $\lambda.$ Its
character is also the Weyl character $s_{\lambda}$. Let $U_{q,\mathbb{Z}}(%
\frak{g})$ be the $\mathbb{Z}[q,q^{-1}]$-subalgebra of $U_{q}(\frak{g})$
generated by the elements 
\begin{equation*}
E_{i}^{(k)}=\frac{E_{i}^{(k)}}{[k]_{i}!},F_{i}^{(k)}=\frac{F_{i}^{(k)}}{%
[k]_{i}!}\text{ and }K_{i}^{\pm1}
\end{equation*}
where $E_{i},F_{i},K_{i}^{\pm1},$ $i\in I_{n}$ are the generators of $U_{q}(%
\frak{g})$. The indeterminate $q$ can be specialized at $\zeta$ in $U_{q,%
\mathbb{Z}}(\frak{g})$. Thus it makes sense to define $U_{\zeta }(\frak{g}%
)=U_{q,\mathbb{Z}}(\frak{g})\otimes_{\mathbb{Z}[q,q^{-1}]}\mathbb{C}$ where $%
\mathbb{Z}[q,q^{-1}]$ acts on $\mathbb{C}$ by $q\mapsto\zeta$. Fix a highest
weight vector $v_{\lambda}$ in $V_{q}(\lambda).$ We have $%
V_{q}(\lambda)=U_{q}(\frak{g})\cdot v_{\lambda}.$ Similarly $%
V_{\zeta}(\lambda)=U_{\zeta}(\frak{g})\cdot v_{\lambda}$ is a $U_{\zeta}(%
\frak{g})$ module called a Weyl module and one has $\mathrm{char}%
(V_{\zeta}(\lambda))=s_{\lambda}.$ The module $V_{\zeta}(\lambda)$ is not
simple but admits a unique simple quotient denoted by $L(\lambda).$

\noindent From results due to Kazhdan-Lusztig and Kashiwara-Tanisaki one
obtains the following decomposition of $\mathrm{char}(L(\lambda))$ on the
basis of Weyl characters:

\begin{theorem}
\label{deo}Consider $\lambda \in P_{+}$.

\begin{enumerate}
\item  For $\ell $ sufficiently large, the character of $L(\lambda )$
decomposes on the form 
\begin{equation}
\mathrm{char}(L(\lambda ))=\sum_{\mu }(-1)^{l(w(\lambda +\rho ))-l(w(\mu
+\rho ))}P_{\mu +\rho ,\lambda +\rho }^{-}(1)s_{\mu }  \label{decKL}
\end{equation}
where the sum runs over the dominant weights $\mu \in P_{+}$ such that $\mu
+\rho \in \widehat{W}\cdot (\lambda +\rho )$.

\item  The parabolic Kazhdan-Lusztig polynomials $P_{\mu +\rho ,\lambda
+\rho }^{-}(q)$ have nonnegative integer coefficients.
\end{enumerate}
\end{theorem}

\noindent \textbf{Remarks:}

\noindent $\mathrm{(i):}$ The decomposition (\ref{decKL}) has been
conjectured by Kazhdan-Lusztig and proved by Kashiwara-Tanisaki.\ In \cite
{KAT}, Kashiwara and Tanisaki have also obtained that the parabolic
Kazhdan-Lusztig polynomials have nonnegative integer coefficients as soon as
the Coxeter system under consideration corresponds to the Weyl group of a
Kac-Moody Lie algebra, thus in the particular context of this paper.

\noindent $\mathrm{(ii):}$ In Assertion 1 of the above theorem, the integer $%
\ell $ is only explictely known in the simply laced case (see Theorem 7.1 in 
\cite{Ta}).

\bigskip

Consider a nonnegative integer $\ell.\;$The Frobenius map $\mathrm{Fr}%
_{\ell} $ is the algebra homomorphism defined from $U_{\zeta}(\frak{g})$ to $%
U(\frak{g})$ by $\mathrm{Fr}_{\ell}(K_{i})=1$ and 
\begin{equation*}
\mathrm{Fr}_{\ell}(E_{i}^{(k)})=\left\{ 
\begin{tabular}{l}
$e_{i}^{(k/\ell)}$ if $\ell$ divides $k$ \\ 
$0$ otherwise
\end{tabular}
\right. \text{ and }\mathrm{Fr}_{\ell}(F_{i}^{(k)})=\left\{ 
\begin{tabular}{l}
$f_{i}^{(k/\ell)}$ if $\ell$ divides $k$ \\ 
$0$ otherwise
\end{tabular}
\right.
\end{equation*}
where $e_{i},f_{i},$ $i\in I_{n}$ are the Chevalley generators of the
enveloping algebra $U(\frak{g})$. This permits to endow each $U(\frak{g})$%
-module $M$ with the structure of a $U_{\zeta}(\frak{g})$-module $M^{\mathrm{%
Fr}_{\ell}}$. Then we have 
\begin{equation*}
\mathrm{char}(M^{\mathrm{Fr}_{\ell}})=\psi_{\ell}(\mathrm{char}(M))
\end{equation*}
in particular for any $\lambda\in P_{+},$ $\mathrm{char}(V(\lambda )^{%
\mathrm{Fr}_{\ell}})=\psi_{\ell}(s_{\lambda}).$

\noindent Each dominant weight $\lambda\in P_{+},$ can be uniquely
decomposed on the form $\lambda=\overset{r}{\lambda}+\ell\overset{q}{\lambda}
$ where $\overset{r}{\lambda},\overset{q}{\lambda}\in P_{+}$ and $\overset{r%
}{\lambda }=(\overset{r}{\lambda}_{1},...,\overset{r}{\lambda}_{n})$
verifies $0\leq\overset{r}{\lambda}_{i+1}-\overset{r}{\lambda}_{i}<\ell$ for
any $i\in I_{n}$.

\begin{theorem}
\label{thLust2}(Lusztig) The simple $U_{\zeta }(\frak{g})$-module $L(\lambda
)$ is isomorphic to the tensor product 
\begin{equation*}
L(\lambda )\simeq L(\overset{r}{\lambda })\otimes V(\overset{q}{\lambda })^{%
\mathrm{Fr}_{\ell }}.
\end{equation*}
\end{theorem}

\noindent By replacing $\lambda$ by $\ell\lambda$ in the previous theorem,
we have $\overset{r}{\lambda}=0$ and $\overset{q}{\lambda}=\lambda.\;$Thus $%
L(\ell\lambda)\simeq V(\lambda)^{\mathrm{Fr}_{\ell}}$. Then one deduces from
(\ref{decKL}) the equality 
\begin{equation*}
\psi_{\ell}(s_{\lambda})=\mathrm{char}(L(\ell\lambda))=\sum_{\mu+\rho \in%
\widehat{W}\cdot(\ell\lambda+\rho)}(-1)^{l(w(\lambda+\rho))-l(w(\mu+\rho
))}P_{\mu+\rho,\ell\lambda+\rho}^{-}(1)s_{\mu}
\end{equation*}
which shows that the coefficients of the expansion of $\psi_{\ell}(s_{%
\lambda })$ on the basis of Weyl characters are, up to a sign, parabolic
Kazhdan-Lusztig polynomials specialized at $q=1.$ This gives 
\begin{equation*}
\left| <\psi_{\ell}(s_{\lambda}),s_{\mu}>\right| =\left| <s_{\lambda
},\varphi(s_{\mu})>\right| =P_{\mu+\rho,\ell\lambda+\rho}^{-}(1).
\end{equation*}
By definition of the action of $\widehat{W}$ on $P$ we have $\widehat{W}%
\cdot(\ell\lambda+\rho)=\widehat{W}\cdot\rho.\;$This implies the

\begin{corollary}
\label{cor}(of Theorems \ref{th_decdeo} and \ref{thLust2}). For $\ell $ as
in 1 of Theorem \ref{th_decdeo}, we have 
\begin{equation*}
\psi _{\ell }(s_{\lambda })=\sum_{\mu +\rho \in \widehat{W}\cdot \rho
}(-1)^{l(w(\lambda +\rho ))-l(w(\mu +\rho ))}P_{\mu +\rho ,\ell \lambda
+\rho }^{-}(1)s_{\mu }.
\end{equation*}
In particular $\varphi (s_{\mu })\neq 0$ if and only if $\mu +\rho \in 
\widehat{W}\cdot \rho ,$ that is $\mu +\rho =w\cdot \rho -\ell \beta $ with $%
w\in W$ and $\beta \in P.$
\end{corollary}

\noindent\textbf{Remark: }The equivalence $\varphi(s_{\mu})\neq
0\Longleftrightarrow\mu+\rho\in\widehat{W}\cdot\rho$ can also be obtained
more elementary from algorithms described in \ref{subsec_comp}.

\subsection{Parabolic K-L polynomials and branching coefficients}

\noindent \textbf{Warning:}\emph{\ In the sequel of the paper we will
suppose }$\ell $\emph{\ sufficiently large so that assertion }1\emph{\ of
Theorem \ref{th_decdeo} holds.\ Moreover }$\ell $\emph{\ is assumed odd when
the Lie groups under consideration are of type }$C$\emph{\ or }$D.$

\bigskip

\noindent Under these hypotheses we have for any $\mu \in \mathcal{P}_{n}$ $%
\varphi _{\ell }(s_{\mu })=0$ or 
\begin{equation}
\varphi _{\ell }(s_{\mu })=\varepsilon (w_{0})S_{\binom{\mu }{\ell },%
\mathcal{I}}  \label{exp_phimu}
\end{equation}
according to the results of \ref{subsec_comp}.

\noindent \textbf{Remark: }According to the algorithms described in \ref
{subsec_comp}, when $\varphi _{\ell }(s_{\mu })\neq 0$, the cardinalities of
the sets $I^{(k)}$ or $X^{(k)}$ contained in $\mathcal{I}$ are determined by
those of the sets $J^{(k)}.\;$In particular they depend only on $n$ and $%
\ell $ and not on the partition $\mu $ considered. Thus in (\ref{exp_phimu}%
), the underlying subgroup of Levi type $G_{\mathcal{I}}$ is, up to
isomorphism, independent on $\mu .$

\bigskip

By using Proposition \ref{prop_dec_Smu} and Theorems \ref{thA}, \ref{thC}, 
\ref{thD}, \ref{thB1}, \ref{thB2} we deduce from Corollary \ref{cor} the

\begin{theorem}
\label{prop_dual}For any $\lambda ,\mu \in \mathcal{P}_{n}$ such that $\mu
+\rho \in \widehat{W}\cdot \rho $%
\begin{equation*}
P_{\mu +\rho ,\ell \lambda +\rho }^{-}(1)=[V(\lambda ):V_{\mathcal{I}}\binom{%
\mu }{\ell }]
\end{equation*}
where $\binom{\mu }{\ell }\;$and $\mathcal{I}$ are obtained from $\mu $ and $%
\ell $ by applying the algorithms described in \ref{subsec_comp}.
\end{theorem}

\subsection{The functions $H_{\protect\mu}^{\ell}$}

For any $\mu\in\mathcal{P}_{n}$, we define the function $G_{\mu}^{\ell}$ by
setting 
\begin{equation}
G_{\mu}^{\ell}=\sum_{\lambda\in\mathcal{P}_{n}}[V(\lambda):V_{\mathcal{I}}%
\binom{\mu}{\ell}]_{q}s_{\lambda}  \label{U}
\end{equation}
where for any $\lambda\in\mathcal{P}_{n}$, $[V(\lambda):V_{\mathcal{I}}%
\binom{\mu}{\ell}]_{q}=P_{\mu+\rho,\ell\lambda+\rho}^{-}(q).$ We also
consider the function $H_{\mu}^{\ell}$ such that 
\begin{equation}
H_{\mu}^{\ell}=G_{\ell\mu}^{\ell}.  \label{H}
\end{equation}

\begin{theorem}
\label{Th_last}Consider a partition $\mu \in \mathcal{P}_{n}.$

\begin{enumerate}
\item  The coefficients of $G_{\mu }^{\ell }$ and $H_{\mu }^{\ell }$ on the
basis of Weyl characters are polynomials in $q$ with nonnegative integer
coefficients.

\item  We have $H_{\mu }^{1}=s_{\mu }$

\item  For $\ell $ sufficiently large $H_{\mu }^{\ell }=Q_{\mu }^{\prime },$
that is $H_{\mu }^{\ell }$ coincide with the Hall-Littlewood function
associated with $\mu .$
\end{enumerate}
\end{theorem}

\noindent To prove our theorem we need the following Lemma:

\begin{lemma}
Consider $\beta \in \mathbb{Z}^{n}.\;$

\begin{itemize}
\item  In type $A_{n-1}$, suppose $\ell >n$.\ Then the weight $\ell \beta
+\rho $ is regular.

\item  In type $B_{n},C_{n}$ or $D_{n},$ suppose $\ell >2n$.\ Then the
weight $\ell \beta +\rho $ is regular.
\end{itemize}
\end{lemma}

\begin{proof}
Consider $w\in W$ and $t_{\gamma }$ such that $t_{\gamma }w\cdot (\ell \beta
+\rho )=\ell \beta +\rho .$ Then $\delta =\ell \beta +\rho -w\cdot (\ell
\beta +\rho )\in (\ell \mathbb{Z)}^{\ell }$. Set $\beta =(\beta
_{1},...,\beta _{n}).\;$For any $i=1,...,n,$ the $i$-th coordinate of $%
\delta $ is $\delta _{i}=\ell \beta _{i}+\rho _{i}-\ell \beta _{w(i)}-\rho
_{w(i)}$. Since $\delta _{i}\in \ell \mathbb{Z}$, we must have $\left| \rho
_{i}-\rho _{w(i)}\right| \in \ell \mathbb{Z}$. One verifies easily that for
type $A_{n-1},$ $\left| \rho _{i}-\rho _{w(i)}\right| \leq n-1$ and for
types $B_{n},C_{n},D_{n}$ $\left| \rho _{i}-\rho _{w(i)}\right| \leq 2n$.
Hence when the conditions of the lemma are verified,$\left| \rho _{i}-\rho
_{w(i)}\right| =0$ for any $i=1,...,n.$ This gives $w=1$ .\ The equality $%
t_{\gamma }w\cdot (\ell \beta +\rho )=\ell \beta +\rho $ implies then that $%
\gamma =0.\;$Thus the stabilizer of $\ell \beta +\rho $ is reduced to $\{1\}$%
, that is $\ell \beta +\rho $ is regular.
\end{proof}

\bigskip

\begin{proof}
(of Theorem \ref{Th_last})

\noindent $1:$ Follows from Theorem \ref{deo} and (\ref{U}).

\noindent $2:$ When $\ell =1$, we have seen that $G=G_{\mathcal{I}}$ and $%
\binom{\mu }{\ell }=\mu .\;$Thus $[V(\lambda ):V_{\mathcal{I}}\binom{\mu }{%
\ell }]_{q}\neq 0$ only if $\lambda =\mu $.\ In this case $H_{\mu
}^{1}=s_{\mu }$ for $[V(\lambda ):V_{\mathcal{I}}\binom{\mu }{\ell }%
]_{q}=[V(\lambda ):V(\lambda )]_{q}=1$.

\noindent $3:$ Suppose $\ell $ as in the previous lemma.\ We have $%
[V(\lambda ):V_{\mathcal{I}}\binom{\ell \mu }{\ell }]_{q}=P_{\ell \mu +\rho
,\ell \lambda +\rho }^{-}(q).$ Since $\ell \lambda +\rho $ is regular for
the action of $\widehat{W}$, we obtain by Theorem \ref{th_decdeo}, $P_{\ell
\mu +\rho ,\ell \lambda +\rho }^{-}(q)=p_{w(\ell \mu +\rho ),w(\ell \lambda
+\rho )}(q)$. By using Lemma \ref{Lem_nlambda}, we deduce $P_{\ell \mu +\rho
,\ell \lambda +\rho }^{-}(q)=p_{n_{\mu },n_{\lambda }}(q).$ Now by Theorem 
\ref{THL} this gives $P_{\ell \mu +\rho ,\ell \lambda +\rho
}^{-}(q)=K_{\lambda ,\mu }(q).$ Finally 
\begin{equation}
H_{\mu }^{\ell }=\sum_{\lambda \in \mathcal{P}_{n}}[V(\lambda ):V_{\mathcal{I%
}}\binom{\ell \mu }{\ell }]_{q}s_{\lambda }=\sum_{\lambda \in \mathcal{P}%
_{n}}K_{\lambda ,\mu }(q)s_{\lambda }=Q_{\mu }^{\prime }.  \label{last}
\end{equation}
\end{proof}

\bigskip

\noindent\textbf{Remarks:}

\noindent$\mathrm{(i):}$ By the previous theorem the functions $%
H_{\mu}^{\ell }$ interpolate between the Weyl characters and the
Hall-Littlewood functions.

\noindent$\mathrm{(ii):}$ When $\ell$ is even for types $C$ and $D$, one can
also define the functions $G_{\mu}^{\ell}$ and $H_{\mu}^{\ell}$ by setting $%
G_{\mu}^{\ell}=\sum_{\lambda\in\mathcal{P}_{n}}P_{\mu+\rho,\ell\lambda+\rho
}^{-}(q)s_{\lambda}$ and $H_{\mu}^{\ell}=G_{\ell\mu}^{\ell}$, respectively.
When $\ell>2n$ we have yet $H_{\mu}^{\ell}=Q_{\mu}^{\prime},$ but the
polynomials $P_{\mu+\rho,\ell\lambda+\rho}^{-}(q)$ cannot be interpreted as
quantizations of branching coefficients.

\noindent$\mathrm{(iii):}$ The conditions $\ell>n$ for type $A_{n-1}$ and $%
\ell>2n$ for types $B_{n},C_{n},D_{n}$ appear also naturally in the
algorithms of \ref{subsec_comp}.\ When they are fulfilled, one has $%
\varphi_{\ell}(s_{\ell\mu})=0,$ or $J_{k}=I_{k}$ for any $k=1,...,n$ and $%
J_{k}=I_{k}=\emptyset$ for $k\notin\{1,...,n\}$. Then $[(\ell\mu)/\ell]=\mu$
and $G_{\mathcal{I}}=H$. Hence $[V(\lambda):V_{\mathcal{I}}\binom{\mu}{\ell }%
]=K_{\lambda,\mu}$ for any $\lambda\in\mathcal{P}_{n}$. This yields equality
(\ref{last}) specialized at $q=1$.

\section{Further remarks\label{Lastsec}}

\subsection{Quantization of tensor product coefficients}

Consider $\mu\in\mathcal{P}_{n}$ and set $\mu=(\mu^{(0)},...,\mu^{(_{%
\ell-1})})$ as in Theorem \ref{thA}.\ For $G=GL_{n},$ the duality $c_{(\mu
^{(0)},...,\mu^{(_{\ell-1})})}^{\lambda}=[V(\lambda):V_{\mathcal{I}}\binom
{\mu}{\ell}]$ yields a $q$-analogue of the Littlewood-Richardson coefficient 
$c_{(\mu^{(0)},...,\mu^{(_{\ell-1})})}^{\lambda}$ defined by setting 
\begin{equation}
c_{(\mu^{(0)},...,\mu^{(_{\ell-1})})}^{\lambda}(q)=[V(\lambda):V_{\mathcal{I}%
}\binom{\mu}{\ell}]_{q}=P_{\mu+\rho,\ell\lambda+\rho}^{-}(q).  \label{quantA}
\end{equation}
By Theorem \ref{deo}, $c_{(\mu^{(0)},...,\mu^{(_{\ell-1})})}^{\lambda}(q)$
have then nonnegative integer coefficients.

\noindent In \cite{lec}, we have shown that there also exists a duality
between tensor product coefficients for types $B,C,D$ defined as the
analogues of the Littlewood-Richardson coefficients by counting the
multiplicities of the isomorphic irreducible components in a tensor product
of irreducible representations and branching coefficients.\ These branching
coefficients correspond to the restriction of $SO_{2n}$ to subgroups of the
form $SO_{2r_{0}}\times \cdot \cdot \cdot SO_{2r_{p}}$ where the $r_{i}$'s
are positive integers summing $n$. These subgroups are not subgroups of Levi
type, thus the Littlewood-Richardson coefficients for types $B,C,D$ cannot
be quantified as in (\ref{quantA}) by using parabolic Kazhdan-Lusztig
polynomials.

\noindent For $G=SO_{2n+1},Sp_{2n}$ or $SO_{2n}$ and $\lambda \in \mathcal{P}%
_{n}$, denote by $\frak{V}(\lambda )$ the restriction of the irreducible
finite dimensional $GL_{N}$-module of highest weight $\lambda $ to $G.$
Consider a $p$-tuple $(\mu ^{(0)},...,\mu ^{(p-1)})$ of partitions such that 
$\mu ^{(k)}\in \mathcal{P}_{r_{k}}$ for any $k=0,...,p-1.\;$One can define
the coefficients $\frak{D}_{\mu ^{(0)},...,\mu ^{(p-1)}}^{\lambda }$ as the
multiplicity of $V(\lambda )$ in $\frak{V}(\mu ^{(0)})\otimes \cdot \cdot
\cdot \otimes \frak{V}(\mu ^{(p-1)}),$ that is such that 
\begin{equation*}
\frak{V}(\mu ^{(0)})\otimes \cdot \cdot \cdot \otimes \frak{V}(\mu
^{(p-1)})\simeq \bigoplus_{\lambda \in \mathcal{P}_{n}}V(\lambda )^{\oplus 
\frak{D}_{\mu ^{(0)},...,\mu ^{(p-1)}}^{\lambda }}.
\end{equation*}
We have also obtained in \cite{lec} a duality result between the
coefficients $\frak{D}_{\mu ^{(0)},...,\mu ^{(p-1)}}^{\lambda }$ and
branching coefficients corresponding to the restriction of $G$ to the
subgroup of Levi type $GL_{r_{0}}\times \cdot \cdot \cdot \times
GL_{r_{p-1}} $. The coefficients $\frak{D}_{\mu ^{(0)},...,\mu
^{(p-1)}}^{\lambda }$ can be expressed by using a partition function
similarly to Proposition \ref{prop_GW}.\ By quantifying this partition
function, one shows that they admit nonnegative $q$-analogues.\ It is
conjectured that stable one-dimensional sums defined in \cite{HKOTY} from
affine crystals obtained by considering the affinizations of the classical
root systems are special cases of the $q$-analogues obtained in this way.

\noindent Recall that the subgroups of Levi type $G_{\mathcal{I}}$ obtained
in the theorems of \ref{subsec_comp} are, up to isomorphism, determined only
by $G$ and $\ell $.\ This implies that there exist subgroups of Levi type $L$
in $G$ which are not isomorphic to a subgroup $G_{\mathcal{I}}$.\ This is
for instance the case when $G=Sp_{2n}$ for the subgroups of Levi type $G_{%
\mathcal{I}}\simeq GL_{r_{0}}\times \cdot \cdot \cdot \times GL_{r_{p-1}}$%
such that $r_{k}>1$ for any $k=0,...,p-1$. Indeed, by Theorem \ref{thC},
when $r_{0}=\mathrm{card}(I^{(0)})>1,$ $G_{\mathcal{I}}$ is isomorphic to 
\begin{equation*}
Sp_{2r_{0}}\times GL_{r_{1}}\times \cdot \cdot \cdot \times GL_{r_{p-1}}.
\end{equation*}
This implies that one cannot obtain in general a quantization of the tensor
product coefficients $\frak{D}_{\mu ^{(0)},...,\mu ^{(p-1)}}^{\lambda }$ by
using parabolic Kazhdan-Lusztig polynomials as in (\ref{quantA}).

\subsection{Combinatorial description of the functions $G_{\protect\mu%
}^{\ell}$}

When $G=GL_{n}$, the functions $G_{\mu}^{\ell}$ defined in (\ref{U}) admit
the following combinatorial description 
\begin{equation*}
G_{\mu}^{\ell}=\sum_{T\in\mathrm{Tab}_{\ell}(\mu)}q^{s(T)}x^{T}
\end{equation*}
where $\mathrm{Tab}_{\ell}(\mu)$ is the set of $\ell$-ribbon tableaux of
shape $\mu$ on $I_{n}$ and $s$ the spin statistic defined on ribbon tableaux%
\textrm{\ }(see \cite{LLT0} page 1057).\ Recently, Haglund, Haiman and Loehr
have obtained the expansion of the Macdonald polynomials in terms of simple
renormalizations of the LLT polynomials $G_{\mu}^{\ell}$.\ This expansion
yields a combinatorial formula for the Macdonald polynomials \cite{HHL}.

\noindent This suggests to investigate the following combinatorial problem:

\begin{problem}
Find a combinatorial description of the polynomials $G_{\mu }^{\ell }$ and
the $q$-analogues $[V(\lambda ):V_{\mathcal{I}}\binom{\mu }{\ell }]_{q}$
related to the roots systems of type $B,C$ or $D$.
\end{problem}

\subsection{Exceptional root systems}

It is also possible to define the plethysm $\psi _{\ell }$ and the dual
plethysm $\varphi _{\ell }$ for exceptional root systems. Consider such an
exceptional root system $R$ and $\mu $ a dominant weight for $R$.$\;$Denote
also by $s_{\mu }$ the Weyl character of the irreducible finite dimensional
module of highest weight $\lambda .\;$When $\ell $ is sufficiently large
(the bound depends on $R$), we have $\varphi _{\ell }(s_{\mu })=s_{\mu }.\;$%
For the other values of $\ell $, one shows that the polynomial $\varphi
_{\ell }(e^{\mu }\prod_{\alpha \in R_{+}}(1-e^{\alpha }))$ do not factorize
in general as a product of factors $(1-x^{\beta })$ where $\beta $ is a
positive root.\ This implies that one cannot define generalized
Hall-Littlewood functions for exceptional types by proceeding as in (\ref{H}%
).

\subsection{Stabilized plethysms\label{subsec-stabipleth}}

When $G=Sp_{2n}$ or $SO_{2n}$ and $\ell $ is even, we have seen that the
combinatorial methods of Section \ref{Sec-Plethy} do not permit to obtain
the coefficients of the expansion of the plethysms $\varphi (s_{\lambda })$
on the basis of the Weyl characters. In \cite{lec2}, we show that this
difficulty can be overcome by considering stabilized power sum plethysms,
i.e.\ by assuming $n\geq \ell \left| \lambda \right| .$ Under this
hypothesis, one can indeed prove that the coefficients in the expansion of $%
\varphi (s_{\lambda })$ coincide for $G=SO_{2n+1},Sp_{2n}$ and $SO_{2n}$. So
it suffices to compute them in type $B_{n}$ for which we have a relevant
combinatorial procedure in both cases $\ell $ even and $\ell $ odd.

\bigskip

\noindent \textbf{Note: }\textit{While revising a previous version of this
work \cite{lec1}, I was informed that Grojnowski and Haiman \cite{HG} also
define, in a paper in preparation, generalized Hall-Littlewood polynomials
for reductive Lie groups. Their polynomials are introduced as formal }$q$%
\textit{-characters depending on a subgroup of Levi type. The coefficients
of the corresponding expansion on the basis of the Weyl characters are also
affine parabolic Kazhdan-Lusztig polynomials. As far as the author can see,
the generalization of the Hall-Littlewood polynomials presented in the
present paper satisfies the general definition given in \cite{HG} (see
Definition 5.12). Nevertheless, our combinatorial results based on the study
of the power sum plethysms on Weyl characters are completely independent of
the approach of Grojnowski and Haiman. It also naturally yields the family
of polynomials }$\{G_{\mu }^{\ell }\mid \ell \in \mathbb{N}\}$ \textit{in
the spirit of the original work by Lascoux, Leclerc and Thibon} \cite{LLT0}.


\bigskip

\begin{thebibliography}{99}
\bibitem{GW}  \textsc{G. Goodman, N. R Wallach, }\textit{Representation
theory and invariants of the classical groups}, Cambridge University Press.

\bibitem{HG}  \textsc{G. Grojnowski, M. Haiman, }\textit{Affine Hecke
algebras and positivity of LLT and Macdonald polynomials, }url:
http://math.berkeley.edu/\symbol{126}mhaiman/.

\bibitem{HHL}  \textsc{J. Haglund, M. Haiman, N. Loehr, }\textit{A
combinatorial formula for Macdonald polynomials, }J.\ Amer.\ Math.\ Soc.\ 
\textbf{18}, n$%
{{}^\circ}%
$ 3, 735-761 (2005).

\bibitem{HKOTY}  \textsc{G. Hatayama, A. Kuniba, M. Okado, T. Takagi, Y.
Yamada,}\textit{\ Remarks on fermionic formula, }in N. Jing and K.\ C.\
Misra, eds.\ Recent Developments in Quantum Affine Algebras and Related
Topics, Contemporary Mathematics \textbf{248}, AMS, Providence, 243-291,
(1999).

\bibitem{KAT}  \textsc{G. Kashiwara, A. Tanisaki, }\textit{Parabolic
Kazhdan-Lusztig polynomials and Schubert varieties, }J.\ Algebra, \textbf{249%
}, 306-325 (2002).

\bibitem{KL}  \textsc{D. Kazhdan,}\textit{\ }\textsc{G. Lusztig,}\textit{\
Representations of Coxeter groups and Hecke algebras, Inventiones \textbf{53}%
, 191-213 (1979).}

\bibitem{KT}  \textsc{K. Koike,}\textit{\ }\textsc{I. Terada,}\textit{\
Young diagrammatic methods for the restriction of representations of complex
classical Lie groups to reductive subgroups of maximal rank, Advances in
Mathematics, 79, 104-135 (1990).}

\bibitem{LLT0}  \textsc{A. Lascoux, B.\ Leclerc, J. Y. Thibon, }\textit{%
Ribbon tableaux, Hall Littelwood functions, quantum affine algebras, }J.\
Math.\ Phys.\ \textbf{38}, 1041-1068 (1996).

\bibitem{LT}  \textsc{B.\ Leclerc, J. Y. Thibon, }\textit{%
Littlewood-Richardson coefficients and Kazhdan-Lusztig polynomials, }Advance
Studies in Pure Mathematics \textbf{28}, Combinatorial Methods in
Representation Theory, 155-220 (2000).

\bibitem{lec}  \textsc{C. Lecouvey, }\textit{Quantization of branching
coefficients for classical Lie groups}, to appear in Journal of Algebra,
arXiv math.RT/\textit{0602089}.

\bibitem{lec1}  \textsc{C. Lecouvey, }\textit{Parabolic Kazhdan-Lusztig
polynomials, plethysm and generalized Hall-Littlewood functions for
classical types}, arXiv math.RT/0607038.

\bibitem{lec2}  \textsc{C. Lecouvey, }\textit{Stabilized plethysms for the
classical Lie groups, }arXiv math.RT/\textit{0703514}.

\bibitem{mac}  \textsc{I-G. Macdonald,} \textit{Symmetric functions and Hall
polynomials}, Second edition, Oxford Mathematical Monograph, Oxford
University Press, New York, (1995).

\bibitem{NR}  \textsc{K. Nelsen, A. Ram,} \textit{Kostka-Foulkes polynomials
and Macdonald spherical functions}, Surveys in Combinatorics 2003, C.
Wensley ed., London Math. Soc. Lect. Notes \textbf{307} , Cambridge
University Press, 325--370 (2003).

\bibitem{Ram}  \textsc{A. Ram, }\textit{Weyl group, symmetric functions and
the representation theory of Lie algebras, }Proceedings of the 4th
conference ``Formal Power Series and Algebraic Combinatorics'', Publ.\ LACIM 
\textbf{11}, 327-342 (1992).

\bibitem{So}  \textsc{W. Soergel, }\textit{Kazhdan-Lusztig polynomials and a
combinatorics for tilting modules, }Represent. Theory\textit{\ }1, 83-114
(1997).

\bibitem{Ta}  \textsc{T. Tanisaki, }\textit{Character Formulas of
Kazhdan-Lusztig type, }Fields institute communications, vol \textbf{40},
261-275 (2004).
\end{thebibliography}
\end{document}